\documentclass[11pt]{amsart}

\usepackage{graphicx,xcolor}
\usepackage{bm}

\newcommand{\be}{\begin{equation}}
\newcommand{\ee}{\end{equation}}
\newcommand{\vct}[1]{\bm{#1}}
\newcommand{\mtx}[1]{\mathsf{#1}}

\newcommand{\lsp}{\vspace{3mm}}

\numberwithin{equation}{section}
\theoremstyle{definition}
\newtheorem{remark}{Remark}
\numberwithin{remark}{section}
\newtheorem{definition}{Definition}
\numberwithin{definition}{section}

\begin{document}

\title[High-order Nystr\"om discretization on curves]%
{High-order accurate methods for Nystr\"om discretization of
integral equations on smooth curves in the plane}

\author{S.~Hao}
\address{Department of Applied Mathematics, University of Colorado at Boulder}
\author{A.~H.~Barnett}
\address{Department of Mathematics, Dartmouth College}
\author{P.~G.~Martinsson}
\address{Department of Applied Mathematics, University of Colorado at Boulder}
\author{P.~Young}

\begin{abstract}
Boundary integral equations and Nystr\"om discretization provide a powerful
tool for the solution of Laplace and Helmholtz boundary value problems.
However, often a weakly-singular kernel arises, in which case
specialized quadratures that modify the matrix entries near the diagonal
are needed to reach a high accuracy.
We describe the construction of four different quadratures
which handle logarithmically-singular kernels.
Only smooth boundaries are considered, but some of the techniques
extend straightforwardly to the case of corners.
Three are modifications of the global periodic trapezoid rule,
due to Kapur--Rokhlin, to Alpert, and to Kress.
The fourth is a modification to a quadrature based on
Gauss-Legendre panels due to Kolm--Rokhlin;
this formulation allows adaptivity. 
We compare in numerical experiments the convergence of the four schemes in
various settings, including low- and high-frequency planar Helmholtz problems,
and
3D axisymmetric Laplace problems.
We also find striking differences in performance in an iterative
setting.
We summarize the relative advantages of the schemes.
\end{abstract}
\maketitle

\section{Introduction} 

Linear elliptic boundary value problems (BVPs)
where the partial differential equation has constant
or piecewise-constant coefficients
arise frequently in engineering, mathematics, and physics.
For the Laplace equation, applications include
electrostatics, heat and fluid flow, and probability;
for the Helmholtz equation they include
the scattering of waves in acoustics,
electromagnetics, optics, and quantum mechanics.
Because the fundamental solution (free-space Green's function) is known,
one may solve such problems
using boundary integral equations (BIEs).
In this approach,
a BVP in two dimensions (2D) is converted via so-called jump relations
to an integral equation for an unknown function
living on a 1D curve \cite{coltonkress}.
The resulting reduced dimensionality and geometric simplicity
allows for high-order accurate numerical solutions
with much more efficiency than standard finite-difference or finite
element discretizations \cite{atkinson1997}.

The BIEs that arise in this setting often take the second-kind form
\begin{equation}
\label{eq:fredholm_bie}
\sigma(x) + \int_{0}^{T} k(x,x')\sigma(x')\,dx' = f(x),\qquad x \in [0,T],
\end{equation}
where $[0,T]$ is an interval, where $f$ is a given smooth $T$-periodic function,
and where $k$ is a (doubly) $T$-periodic kernel function
that is smooth away from the origin and has a logarithmic singularity
as $x' \rightarrow x$.
In order to solve a BIE such as (\ref{eq:fredholm_bie}) numerically,
it must be turned into a linear system with a finite number $N$ unknowns.
This is most easily done via the Nystr\"om method \cite{nystrom,LIE}.
(There do exist other discretization methods such as Galerkin
and collocation \cite{LIE}; while their relative merits
in higher dimensional settings are still debated,
for curves in the plane there seems to be little to compete with
Nystr\"om \cite[Sec.~3.5]{coltonkress}.)
However, since the fundamental solution in 2D has a logarithmic
singularity, generic integral operators of interest inherit this singularity
at the diagonal,
giving them (at most weakly-singular) 
kernels which we write in the standard ``periodized-log'' form
\begin{equation}
\label{eq:assump}
k(x,x') = \varphi(x,x')
\log\left(4\sin^2\frac{\pi(x - x')}{T}\right)
 + \psi(x,x')
\end{equation}
for some smooth, doubly $T$-periodic functions $\varphi$ and $\psi$.
In this note we focus on a variety of high-order quadrature schemes for the
Nystr\"om solution of such 1D integral equations.

We first review the Nystr\"om method, and classify some quadrature schemes
for weakly-singular kernels.


\begin{table}           
\small
\begin{tabular}{l|l|l}
\rule[-2ex]{0ex}{0ex}%
&split into $\varphi$, $\psi$ explicit & split into $\varphi$, $\psi$ unknown\\
\hline
\rule{0ex}{3ex}%
global& $\bullet$ Kress$^\dag$ \cite{kress91} & $\bullet$ Kapur--Rokhlin \cite{Kapur:97a}\\
(periodic trapezoid rule)&&$\bullet$ Alpert \cite{alpert:99a}\\
\rule[-2ex]{0ex}{0ex}%
&&$\circ$ QBX$^\ast$ \cite{qbx}\\
\hline
\rule{0ex}{3ex}%
panel-based& $\circ$ Helsing \cite{helsingmixed,helsingtut} &
$\bullet$ Modified Gaussian (Kolm-Rokhlin) \cite{Kolm:01a}\\
(Gauss-Legendre nodes)&&$\circ$ QBX$^\ast$ \cite{qbx}
\end{tabular}
\vspace{2ex}
\caption{Classification of
Nystr\"om quadrature schemes for logarithmically-singular kernels on
smooth 1D curves.
Schemes tested in this work are marked by a solid bullet (``$\bullet$'').
Schemes are amenable to the FMM unless indicated with a $\dag$.
Finally, $\ast$ indicates that other analytic knowledge is required, namely
a local expansion for the PDE.
\label{t:schemes}
}
\end{table}

\subsection{Overview of Nystr\"{o}m discretization}

\begin{figure}   
(a)\raisebox{-1.4in}{\includegraphics[width=.28\textwidth]{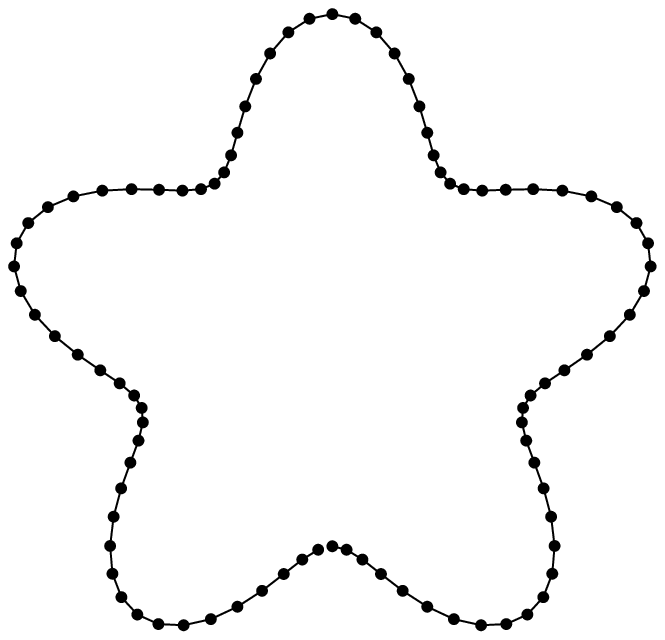}}
(b)\raisebox{-1.4in}{\includegraphics[width=.28\textwidth]{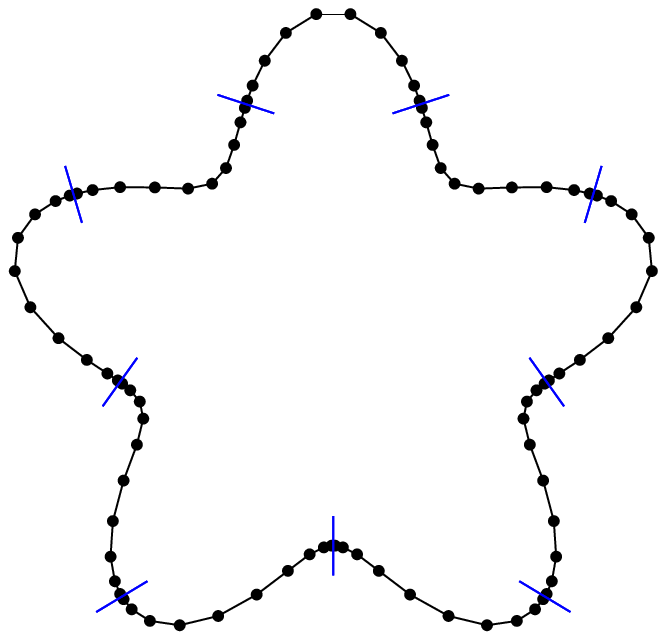}}
\hfill
(c)\raisebox{-1.8in}{\includegraphics[width=.32\textwidth]{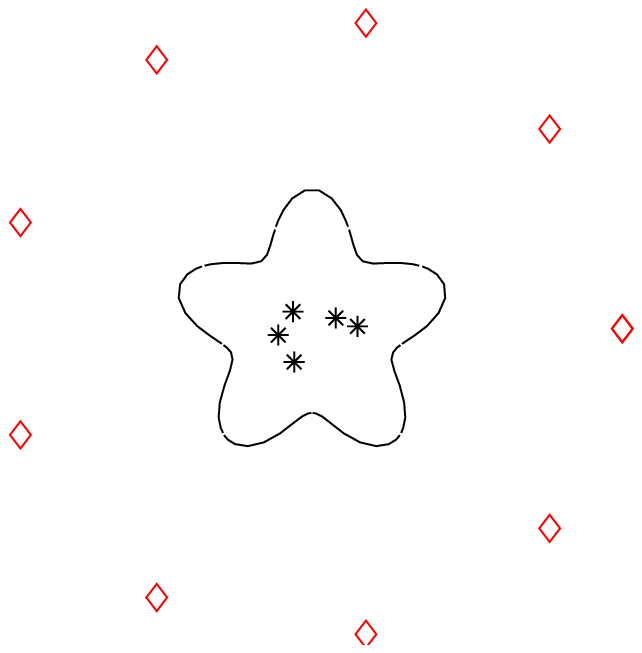}}
\caption{Example smooth planar curve discretized with $N=90$ points via
(a) periodic trapezoid rule nodes and (b) panel-based rule
(10-point Gauss-Legendre; the panel ends are shown by line segments).
In both cases the parametrization is polar angle $t\in[0,2\pi]$
and the curve is the radial function $f(t) = 9/20 - (1/9)\cos(5 t)$.
(c) Geometry for 2D Helmholtz numerical examples in section
\ref{sec:BIE_helm} and \ref{sec:GMRES}. The curve is as in (a) and (b).
Stars show source locations that generate the potential, while
diamonds show testing locations.
\label{fig:2Ddomains}
}
\end{figure}

One starts with an {\em underlying} quadrature scheme
on $[0,T]$, defined by nodes $\{x_{i}\}_{i=1}^{N}$ ordered by
$0 \leq x_{1} < x_{2} < x_{3} < \cdots < x_{N} < T$,
and corresponding weights
$\{w_{i}\}_{i=1}^{N}$.
This means that for $g$ a smooth $T$-periodic function,
$$
\int_0^T g(x) dx \; \approx \; \sum_{i=1}^N w_i g(x_i)
$$
holds to high accuracy. More specifically,
the error converges to zero to high order in $N$.
Such quadratures fall into two popular types:
either a {\em global} rule on $[0,T]$,
such as the periodic trapezoid rule \cite[Sec.~12.1]{LIE}
(which has equally-spaced nodes and equal weights),
or a {\em panel-based} (composite) rule which
is the union of
simple quadrature rules on disjoint intervals (or {\em panels})
which cover $[0,T]$.
An example of the latter is composite Gauss-Legendre quadrature.
The two types are shown in Figure \ref{fig:2Ddomains} (a) and (b).
Global rules may be excellent---for instance,
if $g$ is analytic in a neighborhood of the real axis,
the periodic trapezoid rule has exponential convergence
\cite[Thm.~12.6]{LIE}---%
yet panel-based rules can be more useful in practice because
they are very simple to make adaptive:
one may split a panel into two smaller panels
until a local convergence criterion is met.

The Nystr\"om method for discretizing \eqref{eq:fredholm_bie}
constructs a linear system that
relates a given data vector $\vct{f} = \{f_{i}\}_{i=1}^{N}$
where $f_{i} = f(x_{i})$ to an unknown solution vector $\vct{\sigma} = \{\sigma_{i}\}_{i=1}^{N}$ where
$\sigma_{i} \approx \sigma(x_{i})$.
Informally speaking, the idea is to use the nodes
$\{x_{i}\}_{i=1}^{N}$ as {\em collocation points} where \eqref{eq:fredholm_bie}
is enforced:
\begin{equation}
\label{eq:park1}
\sigma(x_{i}) + \int_{0}^{T} k(x_{i},x')\sigma(x')\,dx' = f(x_{i}),\qquad i = 1,\dots,N.
\end{equation}
Then matrix elements 
$\{a_{i,j}\}_{i,j=1}^{N}$ are constructed such that, for smooth
$T$-periodic $\sigma$,
\begin{equation}
\label{eq:park2}
\int_{0}^{T} k(x_{i},x')\sigma(x')\,dx' \approx \sum_{j=1}^{N} a_{i,j}\,\sigma(x_{j})
~.
\end{equation}
Combining (\ref{eq:park1}) and (\ref{eq:park2}) we obtain a square linear system
that relates $\vct{\sigma}$ to $\vct{f}$:
\begin{equation}
\label{eq:park3}
\sigma_{i} + \sum_{j=1}^{N} a_{i,j}\,\sigma_{j} = f_{i},\qquad i = 1,\dots,N~.
\end{equation}
In a sum such as (\ref{eq:park3}), it is
convenient to think of $x_j$ as the {\em source} node, and $x_i$ as the
{\em target} node.
We write (\ref{eq:park3}) in matrix form as
\begin{equation}
\label{eq:park4}
\vct{\sigma} + \mtx{A}\vct{\sigma} = \vct{f}
~.
\end{equation}
A high-order approximation to $\sigma(x)$ for general $x\in[0,T]$
may then be constructed by interpolation through the values $\vct{\sigma}$.

If the kernel $k$ is smooth, as occurs for the Laplace double-layer operator,
then the matrix elements
\begin{equation}
\label{eq:standquad}
a_{i,j} = k(x_{i},x_{j})w_{j}
\end{equation}
lead to an error convergence rate that is provably the same order
as the underlying quadrature scheme \cite[Sec.~12.2]{LIE}.
It is less obvious how to construct the matrix
$\mtx{A} = \{a_{i,j}\}$ such that (\ref{eq:park2}) holds
to high order accuracy in the case where $k$ has a logarithmic singularity,
as in \eqref{eq:assump}.
The purpose of this note is to describe and compare several
techniques for this latter task.
Note that it is the assumption that the solution $\sigma(x)$ is
smooth (i.e.\ well approximated by high-order interpolation schemes)

\subsection{Types of singular quadrature schemes}

We now overview the schemes presented in this work.
It is desirable for a scheme for the singular kernel case to have {\em almost
all} elements be given by \eqref{eq:standquad}, for the following reason.
When $N$ is large (greater than $10^4$, say),
solving \eqref{eq:park4} via standard dense linear algebra
starts to become impractical, since $O(N^3)$ effort is needed.
Iterative methods
are preferred which converge to a solution using
a small number of matrix-vector products; in the last couple of decades
so-called {\em fast algorithms} have arisen to perform such a
matrix-vector product involving a dense $N\times N$ matrix in only $O(N)$ or $O(N \log N)$ time.
The most well-known is probably the fast multipole
method (FMM) of Rokhlin--Greengard \cite{rokhlin1987}, but others
exist \cite{1986_barnes_hut,2001_krasny_treecode,kifmm}.
They use potential theory to evaluate all $N$ sums of the form
\begin{equation}
\label{eq:FMMsum}
\sum_{j=1}^{N}k(x_{i},x_{j})\,q_{j},\qquad
i = 1,\dots,N
\end{equation}
where the $q_j$ are interpreted as charge strengths.
Choosing $q_j = w_j \sigma_j$ turns this into a fast algorithm to evaluate
$\mtx{A}\vct{\sigma}$ given $\vct{\sigma}$,
in the case where $\mtx{A}$ is a standard Nystr\"om matrix \eqref{eq:standquad}.
\begin{definition}
We say that a quadrature scheme is \textit{FMM-compatible} provided that only
$O(N)$ elements $\{a_{i,j}\}$ differ from the standard formula \eqref{eq:standquad}.
\end{definition}
\noindent
An FMM-compatible scheme can easily be combined with any fast summation scheme
for the sum (\ref{eq:FMMsum}) without compromising its asymptotic speed.
Usually, for FMM-compatible schemes, the elements which differ from \eqref{eq:standquad}
will lie in a band about the diagonal; the width of the band depends only
on the order of the scheme (not on $N$).
All the schemes we discuss are FMM-compatible, apart from that of Kress
(which is not to say that Kress quadrature is incompatible with fast summation;
merely that a standard FMM will not work out of the box).

Another important distinction is the one between (a) schemes in
which the analytic split \eqref{eq:assump} into two smooth kernels must be
{\em explicitly known} (i.e.\ the functions $\varphi$ and $\psi$
are independently evaluated), and (b) schemes which
merely need to access the overall kernel function $k$.
The latter schemes are more flexible, since in applications
the split is not always available (as in the axisymmetric example of
section~\ref{sec:3dBIE}). However, as we will see, this flexibility comes
with a penalty in terms of accuracy.

The following schemes will be described:

\begin{itemize}
\item{\bf Kapur--Rokhlin} (section \ref{sec:KR}).
This is the simplest scheme to implement,
based upon an underlying periodic trapezoid rule.
The weights, but not the node locations, are modified near the diagonal.
No explicit split is needed.

\item {\bf Alpert} (section \ref{sec:alpert}).
Also based upon the periodic trapezoid rule, and also not needing
an explicit split, this scheme replaces
the equi-spaced nodes near the diagonal with an optimal
set of auxiliary nodes, at which new kernel evaluations are needed.

\item {\bf Modified Gaussian} (section \ref{sec:gauss}).
Here the underlying quadrature is Gauss-Legendre panels,
and new kernel evaluations are needed at each set of auxiliary nodes
chosen for each target node in the panel.
These auxiliary nodes are chosen using the algorithm of Kolm--Rokhlin
\cite{Kolm:01a}.
No explicit split is needed.

\item {\bf Kress} (section \ref{sec:kress}).
This scheme uses an explicit split to create a spectrally-accurate
product quadrature based upon the periodic trapezoid rule nodes.
All of the matrix elements differ from the standard form \eqref{eq:standquad},
thus the scheme is not FMM-compatible.
We include it as a benchmark where possible.
\end{itemize}

\noindent
Table~\ref{t:schemes} classifies these schemes (and a couple of others),
according to whether they have underlying global (periodic trapezoid rule)
or panel-based quadrature, whether the split into the
two smooth functions need be explicitly known or not, and
whether they are FMM-compatible.

In section~\ref{sec:num} we present numerical tests comparing the accuracy of
these quadratures in 1D, 2D, and 3D axisymmetric settings.
We also demonstrate that some schemes have negative effects on the
convergence rate in an iterative setting.
We compare the advantages of the
schemes and draw some conclusions in section~\ref{sec:conc}.


\subsection{Related work and schemes not compared}
The methods described in this paper rely on earlier work
\cite{kress91,alpert:99a,Kapur:97a,Kolm:01a}
describing high-order quadrature rules for integrands with
weakly singular kernels. It appears likely that these rules
were designed in part to facilitate Nystr\"om discretization
of BIEs, but, with the exception of Kress \cite{kress91},
the original papers leave most details out.
(Kress describes the Nystr\"om implementation but does not motivate
the quadrature formula; hence we derive this in section~\ref{sec:kress}.)
Some later papers reference the use (e.g.~\cite{2011_bremer,Martinsson:04a})
of high order quadratures but provide few details.
In particular, there appears to have been no formal comparison
between the accuracy of different approaches.

There are several schemes that we do not have the space to include
in our comparison. One of the most promising is the recent scheme of
Helsing for Laplace \cite{helsingmixed} and Helmholtz \cite{helsingtut}
problems,
which is panel-based but uses an explicit split in the style of
Kress, and thus needs no extra kernel evaluations.
We also note the recent QBX scheme \cite{qbx} (quadrature by expansion)
which makes use of off-curve evaluations and local expansions of the PDE.

\section{A brief review of Lagrange interpolation}
\label{sec:math}

This section reviews some well-known (see, e.g., \cite[Sec 3.1]{atkinson1997})
facts about polynomial interpolation that will be used repeatedly in the text.

For a given set of distinct nodes $\{x_{j}\}_{j=1}^N$ and function values $\{y_{j}\}_{j=1}^N$,
the Lagrange interpolation polynomial $L(x)$ is the unique polynomial of degree no greater than
$N-1$  that passes through the $N$ points $\{(x_j, \, y_j)\}_{j=1}^N$. It is given by
$$
L(x) = \sum_{j=1}^N y_j \, L_{j}(x),
$$
where
\begin{equation}
\label{eq:laglag}
L_{j}(x) = \prod_{\substack{i=1 \\ i\neq j}}^N \left( \frac{x - x_{i}}{x_{j} - x_{i}} \right).
\end{equation}
While polynomial interpolation can in general be highly unstable, it is perfectly
safe and accurate as long as the interpolation nodes are chosen well. For instance,
for the case where the nodes $\{x_{j}\}_{j=1}^{N}$ are the nodes associated with
standard Gaussian quadrature on an interval $I = [0,b]$, it is known \cite[Thm.~3.2]{atkinson1997}
that for any $f \in C^{N}(I)$
$$
\left|f(s) - \sum_{j=1}^{N} L_{j}(s)\,f(x_{j})\right| \leq C\,b^{N}\qquad s \in I,
$$
where
$$
C = \left( {\displaystyle\sup_{s \in [0,\,b]}} |f^{(N)}(s)| \right) / N!
$$

\section{Nystr\"om discretization using the Kapur-Rokhlin quadrature rule}
\label{sec:KR}

\subsection{The Kapur--Rokhlin correction to the trapezoid rule}
\label{sec:KR1}

Recall that the standard $N{+}1$-point trapezoid rule that approximates
the integral of a function $g \in C^\infty[0,T]$ is
$h[g(0)/2 + g(h) + \dots + g(T-h) + g(T)/2]$,
where the node spacing is
$$h = \frac {T} {N}
~,$$
and that it has only 2nd-order accuracy \cite[Thm.~12.1]{LIE}.
The idea of Kapur--Rokhlin \cite{Kapur:97a} is to modify this rule
to make it high-order accurate,
by changing a small number of weights near the interval ends, and by adding some
extra equi-spaced
evaluation nodes $x_j = hj$ for $-k\le j< 0$ and $N< j\le N+k$,
for some small integer $k>0$,
i.e.\ nodes lying just {\em beyond} the ends. 
They achieve this goal for functions $g\in C^{\infty}[0,T]$,
but also for the case where one or more endpoint behavior is singular, as in
\begin{equation}
\label{eq:g1}
g(x) = \varphi(x)s(x) + \psi(x)
~,
\end{equation}
where $\varphi(x), \psi(x) \in C^{\infty}[0,T]$ and $s(x)\in C(0,T)$ is a
known function with integrable singularity at zero, or at $T$, or at both
zero and $T$.
Accessing extra nodes outside of the interval
suppresses the rapid growth with order
of the weight magnitudes that plagued previous corrected trapezoid rules.
However, it is then maybe unsurprising that their results need
the additional assumption that $\varphi, \psi \in C^{\infty}[-hk,T+hk]$.

Since we are interested in methods for kernels of the form \eqref{eq:assump},
we specialize to a periodic integrand with the logarithmic singularity at $x=0$ (and therefore also at $x=T$, making both endpoints singular),
\begin{equation}
g(x) = \varphi(x)\log\left|\sin \frac{x\pi}{T}\right| + \psi(x)
~.
\label{glog}
\end{equation}
The $m$th-order Kapur--Rokhlin rule $T_m^{N+1}$ which corrects for
a log singularity at both left and right endpoints is,
\begin{equation}
\label{eq:KRquad}
T_m^{N+1}(g) = h\biggl[\sum_{\substack{\ell=-m \\ \ell \neq 0}}^m \gamma_\ell \, g(\ell h)
+ g(h) + g(2h) + \dots + g(T-h)
+ \sum_{\substack{\ell=-m \\ \ell \neq 0}}^m \gamma_{-\ell} \, g(T + \ell h) \biggr]
\end{equation}
Notive that the left endpoint correction weights $\{\gamma_\ell\}_{\ell=-m, \ell \neq 0}^m$ are also used (in reverse order) at the right endpoint.
The convergence theorems from \cite{Kapur:97a} then imply that,
for any fixed $T$-periodic $\varphi, \psi \in C^{\infty}(\mathbb{R})$,
\begin{equation}
\left|\int_0^T g(x)\,dx - T_m^{N+1}(g)\right| =
O(h^m)\qquad \mbox{ as } N\to\infty
~.
\end{equation}

The apparent drawback that function values are needed at
nodes $\ell h$ for $-m\le \ell\le -1$, which lie
outside the integration interval, is not a problem
in our case of periodic functions, since
by periodicity these function values are known.
This enables us to rewrite \eqref{eq:KRquad} as
\begin{equation}
T_m^{N+1}(g) = h \biggl[\sum_{\ell=1}^m (\gamma_\ell +\gamma_{-\ell}) g(\ell h)
+ g(h) + g(2h) + \dots + g(T-h)
+ \sum_{\ell=-m}^{-1} (\gamma_\ell +\gamma_{-\ell}) g(T + \ell h) \biggr]
\label{eq:KRwrap}
\end{equation}
which involves only the $N-1$ nodes interior to $[0,T]$.

For orders $m=2, 6, 10$
the values of $\gamma_\ell$ are given in the
left-hand column of \cite[Table~6]{Kapur:97a}.
In our periodic case, only the values $\gamma_\ell+\gamma_{-\ell}$ are needed;
for convenience we give them in appendix~\ref{a:KRrules}.
Notice that they are of size 2 for $m=2$, of size around 20 for $m=6$,
and of size around 400 for $m=10$, and alternate in sign in each case.
This highlights the
statement of Kapur--Rokhlin that they were only partially
able to suppress the growth of weights in the case of singular endpoints
\cite{Kapur:97a}.


\subsection{A Nystr\"om scheme}
\label{sec:KR2}
We now can
construct numbers $a_{i, j}$ such that (\ref{eq:park2}) holds.
We start with the underlying periodic trapezoid rule,
with equi-spaced nodes $\{x_{i}\}_{i = 1}^N$ with $x_j = hj$, $h = T/N$.
We introduce a discrete offset
function $\ell(i,j)$ between two nodes $x_{i}$ and $x_{j}$ defined by
$$\ell(i,j) \; \equiv \; j-i \; (\text{mod} \ N),
\qquad -N/2 < \ell(i,j) \le N/2~,
$$
and note that for each $i\in \{1, \ldots, N\}$, and each $|\ell|<N/2$,
there is a unique $j \in \{1,\ldots,N\}$ such that $\ell(i,j) = \ell$.
Two nodes $x_{i}$ and $x_{j}$
are said to be ``close'' if $|\ell(i,j)|\leq m$. Moreover, we call
$x_{i}$ and $x_{j}$ ``well-separated'' if they are not close.

We now apply \eqref{eq:KRwrap} to the integral
(\ref{eq:park2}), which by periodicity of $\sigma$ and the kernel,
we may rewrite with
limits $x_i$ to $x_i + T$ so that the log singularity appears at the
endpoints.
We also extend the definition of the nodes $x_j = hj$ for all integers $j$,
and get,
\begin{align}
\int_0^T k(x_{i}, x')\sigma(x')\,dx' & =
\int_{x_{i}}^{x_{i}+T}k(x_i, x')\sigma(x')\,dx'
\nonumber \\
&\approx
h\sum_{j=i+1}^{i+N-1} k(x_{i}, x_{j})\sigma(x_{j})
    + h\sum_{\substack{\ell=-m\\ \ell\neq 0}}^m (\gamma_{\ell}+\gamma_{-\ell}) k(x_{i}, x_{i+\ell})\,\sigma(x_{i+\ell})~.
\nonumber
\end{align}
Wrapping indices back into $\{1,\ldots,N\}$, the entries of the coefficient matrix $\mtx{A}$ are seen to be,
\begin{equation}
\label{eq:KRentry1}
a_{i,j} = \left\{\begin{array}{ll}
0 \quad&\text{if $i=j$},\\，
h\,k(x_{i}, x_{j}) \quad&\text{if  $x_{i}$ and $x_{j}$ are \lq\lq{}well-separated\rq\rq{}},\\
h\,(1+\gamma_{\ell(i,j)} + \gamma_{-\ell(i,j)})\,k(x_{i}, x_{j})\quad&\text{if $x_{i}$ and $x_{j}$ are \lq\lq{}close\rq\rq{}, and $i\neq j$.}\\
\end{array}\right.
\end{equation}
Notice that
this is the elementwise product of a circulant matrix with the kernel matrix
$k(x_i,x_j)$, and that diagonal values are ignored.
Only $O(N)$ elements (those closest to the diagonal) differ from the standard
Nystr\"om formula \eqref{eq:standquad}.

\section{Nystr\"om discretization using the Alpert quadrature rule}
\label{sec:alpert}

\subsection{The Alpert correction to the trapezoid rule}

Alpert quadrature is another correction to the trapezoid rule
that is high-order accurate for integrands of the
form \eqref{eq:g1} on $(0,T)$.
The main difference with Kapur--Rokhlin is that Alpert quadrature
uses node locations {\em off} the equi-spaced grid $x_j=hj$, but within the interval $(0,T)$.
Specifically, for the function $g$ in \eqref{glog} with log singularities at both ends,
we denote by $S_{l}^{N+1}(g)$ the $l$th-order Alpert quadrature rule
based on an $N{+}1$-point trapezoid grid, defined by the formula
\begin{equation}
\label{eq:alpert_quad}
S_{l}^{N+1}(g) = h\sum_{p = 1}^m w_{p}\,g(\chi_{p}\,h) + h\sum_{j=a}^{N-a}g(jh)
                         + h\sum_{p = 1}^m w_{p}\,g(T-\chi_{p}\,h).
\end{equation}
There are $N-2a+1$ internal equi-spaced nodes with spacing $h = T/N$
and equal weights; these are common to the trapezoid rule.
There are also $m$ new ``correction''
nodes at each end which replace the $a$ original
nodes at each end in the trapezoid rule.
The label ``$l$th-order'' is actually slightly too strong:
the scheme is proven \cite[Cor.~3.8]{alpert:99a}
to have error convergence of order $O(h^l |\log h|)$ as $h\to0$.
The number $m$ of new nodes needed per end is either $l{-}1$ or $l$.
For each order $l$, the integer $a$ is chosen by experiments to be
the smallest integer leading to positive correction nodes and weights.
The following table shows the values of $m$ and $a$ for log-singular kernels for
convergence orders $l = 2,\,6,\,10,\,16$:

\vspace{1mm}

\begin{center}
\begin{tabular}{l|r|r|r|r}
Convergence order                   & $h^{2}|\log h|$ & $h^{6}|\log h|$ & $h^{10}|\log h|$ & $h^{16}|\log h|$ \\ \hline
Number of correction points $m$     & 1               & 5               & 10               & 15 \\
Width of correction window $a$      & 1               & 3               &  6               & 10
\end{tabular}
\end{center}

\vspace{1mm}

\begin{figure}
\includegraphics[width=.65\textwidth, height = 0.08\textwidth]{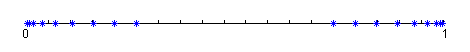}
\caption{Example of Alpert quadrature scheme of order $l=10$
on the interval $[0,1]$.
The original trapezoid rule had 20 points including both endpoints,
i.e.\ $N=19$ and $h=1/19$.
Correction nodes $\{\chi_ph\}_{p=1}^m$ and $\{1-\chi_ph\}_{p=1}^m$
for $m=10$ and $a=6$, are denoted by stars.}
\label{fig:chi}
\end{figure}

\noindent
The corresponding node locations $\chi_{1},\,\dots\,\chi_{m}$ and weights $w_{1},\,\dots\,w_{m}$
are listed in Appendix \ref{a:Alprules}; and illustrated for the case $l=10$ in Figure \ref{fig:chi}.
Details on how to construct these numbers by solving a nonlinear system
can be found in \cite[Sec.~5]{alpert:99a}, and the quadratures listed are
adapted from \cite[Table~8]{alpert:99a}.


\subsection{A Nystr\"om scheme}
Recall that we wish to construct a matrix $a_{i,j}$ that when
applied to the vector of values $\{\sigma(x_j)\}_{j=1}^N$
approximates the action of the integral operator on the function $\sigma$,
evaluated at each target node $x_i$.
To do this via the $l$th-order Alpert scheme with parameter $a$,
we start by using periodicity to shift
the domain of the integral in (\ref{eq:park2})
to $(x_i,x_i+T)$, as in section \eqref{sec:KR2}.
Since the $2m$ auxiliary Alpert nodes lie symmetrically about the singularity
location $x_i$, for simplicity we will treat them as a single set
by defining $\chi_{p+m} = -\chi_p$ and $w_{p+m} = w_p$ for $p=1,\ldots,m$.
Then the rule \eqref{eq:alpert_quad} gives
\be
\label{intAlp}
\int_0^T k(x_{i}, x')\sigma(x')\,dx'
\;\approx \; h\,\sum_{p = a}^{N - a} k(x_i, x_i+ph)\,\sigma(x_i+p h)
+ h\,\sum_{p=1}^{2m}w_{p}k(x_i, x_i+\chi_ph)\,\sigma(x_i+\chi_ph)
~.
\ee

The values $\{\chi_{p}\}_{p=1}^{2m}$ are not integers,
so no auxiliary nodes coincide with any equispaced nodes $\{x_j\}_{j=1}^N$
at which the vector of $\sigma$ values is given.
Hence we must interpolate $\sigma$ to the auxiliary
nodes $\{x_{i}+\chi_{p}h\}_{p=1}^{2m}$.
We do this using local Lagrange interpolation through
$M$ equispaced nodes surrounding the auxiliary source
point $x_{i}+\chi_p h$.
For $M>l$ the
interpolation error is higher order than that of the Alpert
scheme;
we actually use $M=l+3$ since the error is then negligible.
\begin{remark}
While high-order Lagrange interpolation through equi-spaced points is
generally a bad idea due to the Runge phenomenon \cite{tref},
here we will be careful to ensure that the evaluation point always
lies nearly at the center of the interpolation grid, and there is no
stability problem.
\end{remark}
For each auxiliary node $p$,
let the $n$th interpolation node index offset relative to $i$
be
$$
o^{(p)}_n := \lfloor \chi_p - M/2\rfloor + n
~,
$$
and let the whole set be $O^{(p)} := \{o^{(p)}_n\}_{n=1}^M$.
For $q\in O^{(p)}$,
let the function $n_p(q):= q - \lfloor \chi_p - M/2\rfloor$
return the node number of an index offset of $q$.
Finally, let
$$
L^{(p)}_n(x) = \prod_{\substack{k=1 \\ k\neq n}}^M \left(\frac{x-o^{(p)}_k}{o^{(p)}_n - o^{(p)}_k}\right)
$$
be the $n$th Lagrange basis polynomial for auxiliary node $p$.
Applying this interpolation in $\sigma$ gives
for the second term in \eqref{intAlp},
$$
h\,\sum_{p=1}^{2m}w_{p}k(x_i, x_i+\chi_ph)\,\sum_{n=1}^M L_n^{(p)}(\chi_p)
\sigma(x_{i+o_n^{(p)}})
~.
$$
Note that all node indices in the above will be
cyclically folded back into the set $\{1,\ldots, N\}$.

Recalling the notation $\ell(i,j)$ from section~\ref{sec:KR2},
we now find the coefficient matrix $\mtx{A}$ has entries
\begin{equation}
\label{eq:alpert_entry}
a_{i,j} = b_{i,j}+c_{i,j},
\end{equation}
where the first term in \eqref{intAlp} gives
\begin{equation}
\label{eq:b_entry}
b_{i,j} = \left\{\begin{array}{ll}
0 \quad &\text{if $|\ell(i,j)| < a$},\\
h\,k(x_{i}, x_{j}) \quad &\text{if $|\ell(i,j)| \geq a$},\\
\end{array}\right.
\end{equation}
the  standard Nystr\"om matrix \eqref{eq:standquad} with
a diagonal band set to zero, and the auxiliary nodes give
\begin{equation}
\label{eq:c_entry}
c_{i,j} \;=\; h \!\!\! \sum_{\substack{p=1\\ O^{(p)}\ni \ell(i,j)}}^{2m} \!\!\! w_p \,
k(x_i, x_i+\chi_ph) \,
L^{(p)}_{n_p(\ell(i,j))}(\chi_p)
~.
\end{equation}
Notice that the bandwidth of matrix $c_{i,j}$ does not exceed
$a+M/2$, which is roughly $l$,
and thus only $O(N)$ elements of $a_{i,j}$ differ from
those of the standard Nystr\"om matrix.

\section{Nystr\"om discretization using modified Gaussian quadrature}
\label{sec:gauss}

We now turn to a scheme with panel-based underlying quadrature,
namely the composite Gauss-Legendre rule.
We recall that for any interval $I = [0,b]$, the single-panel $n$-point
Gauss-Legendre rule has nodes $\{x_j\}_{j=1}^n \subset I$ and weights
$\{w_j\}_{j=1}^n \subset (0,\infty)$ such that the identity
\begin{equation}
\int_{0}^{b} f(x)\,dx = \sum_{j=1}^n w_{j}f(x_j)
\end{equation}
holds for every polynomial $f$ of degree at most $2n-1$.
For analytic $f$ the rule is exponentially convergent in $n$ with rate given
by the largest ellipse with foci $0$ and $b$ in which $f$ is analytic
\cite[Ch.~19]{ATAP}.

\subsection{Modified Gaussian quadratures of Kolm--Rokhlin}
\label{sec:gauss1}
Suppose that given an  interval $[0,b]$ and a target point $t\in[0,b]$,
we seek to approximate integrals of the form
\begin{equation}
\label{eq:kForm}
\int_0^b \bigl(\varphi(s) \,S(t,s)+ \psi(s)\bigr)\,ds,
\end{equation}
where $\varphi$ and $\psi$ are smooth functions over $[0, b]$ and $S(t,s)$ has a known singularity as $s\rightarrow t$.
Since the integrand is non-smooth, standard Gauss-Legendre quadrature would be inaccurate
if applied to evaluate (\ref{eq:kForm}). Instead, we seek an $m$-node ``modified
Gaussian'' quadrature rule with weights $\{v_{k}\}_{k=1}^m$ and nodes $\{y_{k}\}_{k=1}^m\subset [0,b]$ that evaluates
the integral (\ref{eq:kForm}) to high accuracy. In particular, we use a quadrature
formula of the form
\begin{equation}
\label{eq:nearquad}
\int_0^b \bigl(\varphi(s) \,S(t,s)+ \psi(s)\bigr)\,ds \approx \sum_{k=1}^m v_{k} \bigl(\varphi(y_{k}) \,S(t,y_{k})+ \psi(y_{k})\bigr)
\end{equation}
which holds when $\varphi$ and $\psi$ are polynomials of degree $n$.
It is crucial to note that $\{v_{k}\}_{k=1}^m$ and $\{y_{k}\}_{k=1}^m$
depend on the target location $t$; unless $t$ values are very close then
new sets are needed for {\em each} different $t$.

Next consider the problem of evaluating (\ref{eq:kForm}) in the case where
the target point $t$ is near the region of integration but not actually inside it,
for instance $t \in [-b,0) \cup (b,\,2b]$. In this case, the integrand is smooth
but has a nearby singularity:
this reduces the convergence rate and
means that its practical accuracy would be low for the fixed
$n$ value we prefer.
In this case, we can approximate the integral (\ref{eq:kForm}) with another set of
modified Gaussian quadrature weights $\{\hat{v}_{k}\}_{k=1}^{m'}$ and nodes
$\{\hat{y}_{k}\}_{k=1}^{m'}\subset [0,b]$ giving a quadrature formula analogous
to (\ref{eq:nearquad}).
In fact, such weights and nodes can be found that hold to high accuracy
for all targets $t$ in intervals of the form $[-10^{-p+1},-10^{-p}]$.

Techniques for constructing such modified quadratures
based upon nonlinear optimization are presented by Kolm--Rokhlin
\cite{Kolm:01a}. In the numerical experiments in section \ref{sec:num},
we use a rule with $n=10$, $m=20$, and $m'=24$. This rule leads to
fairly high accuracy, but is not of so high order that clustering 
of the end points becomes an issue in double precision arithmetic.
Since the full tables of quadrature weights get long in this case,
we provide them as text files at \cite{2012_nystrom_website}.


\subsection{A Nystr\"om scheme}
\label{sec:gauss2}

We partition the domain of integration as
$$
[0, T] = \bigcup_{p=1}^{N_{P}}\Omega_p,
$$
where the $\Omega_p$'s are non-overlapping subintervals called panels.
For simplicity, we for now assume that the panels are equi-sized so that
$\Omega_{p} = \left[\tfrac{T(p-1)}{N_{P}},\,\tfrac{Tp}{N_{P}}\right]$.
Note that an adaptive version would have variable-sized panels.
On each panel, we place the nodes of an $n$-point Gaussian quadrature to obtain
a total of $N = n\,N_{P}$ nodes.
Let $\{x_{i}\}_{i=1}^{N}$ and $\{w_{i}\}_{i=1}^{N}$ denote the nodes and weights
of the resulting composite Gaussian rule.

Now consider the task of approximating the integral \eqref{eq:park2} at a target
point $x_{i}$. We decompose the integral as
\begin{equation}
\int_{0}^{T} k(x_{i},x')\,\sigma(x')\,dx' =
\sum_{q=1}^{N_P}\int_{\Omega_q} k(x_{i},x')\,\sigma(x')\,dx'.
\end{equation}
We will construct an approximation for each panel-wise integral
\begin{equation}
\label{eq:inteomega}
\int_{\Omega_q} k(x_{i},x')\sigma(x')\,dx',
\end{equation}
expressed in terms of the values of $\sigma$ at the Gaussian nodes in the source panel
$\Omega_{q}$. There are three distinct cases:

\lsp

\noindent
\textit{Case 1: $x_i$ belongs to the source panel $\Omega_{q}$:}
The integrand is now singular in the domain of integration, but we can exploit
that $\sigma$ is still smooth and can be approximated via polynomial interpolation on this single panel.
Using the set of $n$ interpolation nodes $\{x_k:x_k\in \Omega_q\}$, let
$L_j$ be the Lagrange basis function corresponding to node $j$.
Then,
\begin{equation}
\label{eq:logitech2}
\sigma(x') \; \approx \sum_{j\;:\;x_{j} \in \Omega_{q}} L_{j}(x')\,\sigma(x_{j})
~.
\end{equation}
Inserting (\ref{eq:logitech2}) into the integral in (\ref{eq:inteomega}) we find
$$
\int_{\Omega_q} k(x_{i},x')\sigma(x')\,dx' \; \approx
\sum_{j\;:\;x_{j} \in \Omega_{q}}
\left(\int_{\Omega_q} k(x_{i},x')  L_{j}(x')\,dx'\right)\,\sigma(x_{j}).
$$
Let $\{v_{i,k}\}_{k=1}^m$ and $\{y_{i,k}\}_{k=1}^m$ be the modified Gaussian weights
and nodes in the rule (\ref{eq:nearquad})
on the interval $\Omega_q$ associated with the target $t=x_i$.
Using this rule,
\begin{equation}
\label{eq:formula_diag}
\int_{\Omega_q} k(x_{i},x')\sigma(x')\,dx' \; \approx
\sum_{j \;:\; x_{j} \in \Omega_{q}}\left(\sum_{k=1}^m v_{i,k}\,k(x_{i}, y_{i,k})\,L_j(y_{i,k})\right)\,\sigma(x_j).
\end{equation}
Note that the expression in brackets
gives the matrix element $a_{i,j}$.
Here the auxiliary nodes $y_{i,k}$ play a similar role to the
auxiliary Alpert nodes $x_i + \chi_ph$ from section~\ref{sec:alpert}:
new kernel evaluations are needed at each of these nodes.

\lsp

\noindent
\textit{Case 2: $x_i$ belongs to a panel $\Omega_{p}$ adjacent to source panel $\Omega_{q}$:}
In this case, the kernel $k$ is smooth, but has a singularity closer to
$\Omega_q$ than
the size of one panel, so standard Gaussian quadrature would still be inaccurate. We therefore proceed
as in Case 1: we replace $\sigma$ by its polynomial interpolant and then integrate using
the modified quadratures described in Section \ref{sec:gauss1}. The end result is a
formula similar to (\ref{eq:formula_diag}) but with the sum including $m'$ rather than $m$
terms, and with $\hat{v}_{i,k}$ and $\hat{y}_{i,k}$ replacing $v_{i,k}$ and $y_{i,k}$, respectively.

\lsp

\noindent
\textit{Case 3: $x_{i}$ is well-separated from the source panel $\Omega_{q}$:}
By ``well-separated'' we mean that $x_i$ and $\Omega_q$ are at least
one panel size apart in the parameter $x$.
(Note that if the curve geometry involves close-to-touching parts, then
this might not be a sufficient criterion for being
well-separated in $\mathbb{R}^2$; in practice this would best be
handled by adaptivity.)
In this case, both the kernel $k$ and the potential $\sigma$ are smooth, so the
original Gaussian rule will be accurate,
\begin{equation}
\label{eq:formula_offd}
\int_{\Omega_q}k(x_i,x')\sigma(x')\,dx' \approx
\sum_{j\;:\;x_{j} \in \Omega_{q}}\,w_{j}k(x_i,x_j)\sigma(x_j)
~.
\end{equation}

\lsp

Combining (\ref{eq:formula_offd}) and (\ref{eq:formula_diag}), we find that
the Nystr\"om matrix elements $a_{i, j}$ are given by
\begin{equation}
\label{eq:aij_highorder}
a_{i,j} = \left\{\begin{array}{ll}
\sum_{k=1}^{m} v_{i,k}\,k(x_{i},y_{i,k})\, L_j(y_{i,k}), \quad&\mbox{if }x_{i}\mbox{ and }x_{j}\mbox{ are in the same panel,}
\\
\sum_{k=1}^{m'} \hat{v}_{i,k}\,k(x_{i},\hat{y}_{i,k})\, L_j(\hat{y}_{i,k}), \quad&\mbox{if }x_{i}\mbox{ and }x_{j}\mbox{ are in adjacent panels,}
\\
\rule{0ex}{2.3ex}  
k(x_{i},x_{j})\,w_{j}, \quad&\mbox{if }x_{i}\mbox{ and }x_{j}\mbox{ are in well-separated panels.}
\end{array}\right.
\nonumber
\end{equation}

\section{Nystr\"om discretization using the Kress quadrature rule}
\label{sec:kress}

The final scheme that we present
returns to an underlying periodic trapezoid rule,
but demands separate knowledge of the smooth kernel functions
$\varphi$ and $\psi$ appearing in \eqref{eq:assump}.
We first review spectrally-accurate product quadratures,
which is an old idea, but which we do not find well explained in the
standard literature.

\subsection{Product quadratures}
For simplicity, and
to match the notation of \cite{kress91}, we fix the period $T=2\pi$
and take $N$ to be even. The nodes are thus
$x_i = 2\pi i/N$, $i=1,\ldots,N$.

A {\em product quadrature} approximates the
integral of the product of a general smooth $2\pi$-periodic real function $f$
with a fixed known (and possibly singular) $2\pi$-periodic real function $g$,
by a periodic trapezoid rule with modified weights $w_j$,
\be
\int_0^{2\pi} f(s) g(s) \, ds \; \approx \; \sum_{j=1}^N w_j f(x_j)
~.
\label{pq}
\ee
Using the Fourier series $f(s) = \sum_{n\in\mathbb{Z}} f_n e^{ins}$, and similar
for $g$, we recognize the integral as an inner product and use Parseval,
\be
\int_0^{2\pi} f(s) g(s) \, ds \; = \; 2\pi \sum_{n\in\mathbb{Z}} f_n \overline{g_n}
\label{fg}
\ee
Since $f$ is smooth, $|f_n|$ decays to zero with a high order as $|n|\to\infty$.
Thus we can make two approximations.
Firstly, we truncate the infinite sum to $\sum'_{|n|\le N/2}$,
where the prime indicates that the extreme terms $n=\pm N/2$ are given a factor of $1/2$.
Secondly, we use the periodic trapezoid rule to evaluate the Fourier
coefficents of $f$, i.e.\
\be
f_n \; = \; \frac{1}{2\pi}\int_0^{2\pi} e^{-ins} f(s) ds \; \approx\;
\frac{1}{N} \sum_{j=1}^N e^{-inx_j} f(x_j)
~.
\label{fn}
\ee
Although the latter introduces aliasing (one may check
that the latter sum is exactly $f_n + f_{n+N} + f_{n-N} + f_{n+2N} + f_{n-2N} + \cdots$),
the decay of $f_n$ means that errors decay to high order with $N$.
Substituting \eqref{fn} into the truncated version of \eqref{fg} gives
\be
\int_0^{2\pi} f(s) g(s) \, ds \; \approx \;
2\pi \sum'_{|n|\le N/2} \overline{g_n}
\frac{1}{N} \sum_{j=1}^N e^{-inx_j} f(x_j)
\; \approx \;
\sum_{j=1}^N \biggl( \frac{2\pi}{N}\sum'_{|n|\le N/2} e^{-inx_j} \overline{g_n} \biggr) f(x_j)
\ee
The bracketed expression gives the weights
in \eqref{pq}.
Since $g$ is real (hence $g_{-n}=\overline{g_n}$),
\be
w_j \;=\;  \frac{2\pi}{N}\!\!\sum'_{|n|\le N/2} \!\!e^{-inx_j} \overline{g_n} \;=\;
\frac{2\pi}{N} \biggl[g_0 + \sum_{n=1}^{N/2-1} 2\mbox{Re}(g_n e^{i n x_j})
+ \mbox{Re}(g_{N/2} e^{i N x_j/2})
\biggr],
\quad j=1,\dots,N.
\label{wj}
\ee

\subsection{The Kress quadrature}
To derive the scheme of Kress (originally due to Martensen--Kussmaul; see
references in \cite{kress91}) we note the Fourier series
(proved in \cite[Thm.~ 8.21]{LIE}),
\be
g(s) = \log\left( 4 \sin^2 \frac{s}{2}\right)
\qquad \Leftrightarrow \qquad
g_n =
\left\{\begin{array}{ll}0,& n=0,\\ -1/|n|, &n \neq 0.\end{array}\right.
\label{logfou}
\ee
Translating $g$ by a displacement $t\in\mathbb{R}$
corresponds to multiplication of $g_n$ by $e^{-int}$.
Substituting this displaced series into \eqref{wj} and simplifying
gives
\begin{equation}
\label{eq:pqlog}
\int_{0}^{2\pi} \log \left(4\sin^2 \frac {t-s}{2}\right)\,\varphi(s)\,ds \; \approx \; \sum_{j=1}^{N} R_{j}^{(N/2)}(t)\,\varphi(x_{j})
~,
\end{equation}
where the weights, which depend on the target location $t$,
are
\be
R_{j}^{(N/2)}(t) \; = \; -\frac{4\pi}{N} \biggl[
\sum_{n=1}^{N/2-1} \frac{1}{n} \cos n(x_j - t)
\; + \; \frac{1}{N}\cos \frac{N}{2}(x_j-t)
\biggr]
~, \quad j=1,\dots,N.
\label{Rj}
\ee
This matches \cite[(3.1)]{kress91} (note the convention
on the number of nodes differs by a factor 2).

When a smooth function is also present, we use the periodic trapezoid
rule for it, to get
\begin{equation}
\label{eq:kressquad}
\int_{0}^{2\pi} \log \left(4\sin^2 \frac {t-s}{2}\right)\,\varphi(s) + \psi(s)
\,ds \; \approx \; \sum_{j=1}^{N} R_{j}^{(N/2)}(t)\,\varphi(x_{j})
+ \frac{2\pi}{N} \sum_{j=1}^{N} \psi(x_j)
~.
\end{equation}
Assuming the separation into $\varphi$ and $\psi$ is known,
this gives a high-order accurate quadrature;
in fact for $\varphi$ and $\psi$ analytic, it is exponentially convergent
\cite{kress91}.

\subsection{A Nystr\"om scheme}
We use the above Kress quadrature to approximate the integral (\ref{eq:park2})
where the kernel has the form \eqref{eq:assump} with $T=2\pi$,
and the functions $\varphi(x,x')$ and $\psi(x,x')$ are separately known.
Applying \eqref{eq:kressquad}, with $h=2\pi/N$, gives
\begin{align}
\int_{0}^{2\pi}k(x_{i}, x')\sigma(x')\,dx' \;\approx\;
\sum_{j=1}^{N}R_{j}^{(N/2)}(x_{i})\varphi(x_{i}, x_{j})\sigma(x_{j}) +
h\sum_{j=1}^{N}\psi(x_{i}, x_{j})\sigma(x_{j})
~.
\end{align}
Using the symbol $R_{j}^{(N/2)} := R_{j}^{(N/2)}(0)$,
and noticing that $R_{j}^{(N/2)}(x_{i})$ depends only on $|i-j|$, we
find that the entries of the coefficient matrix $\mtx{A}$ are,
\begin{equation}
\label{eq:kressentry}
a_{i,j} \; =\;
R_{|i-j|}^{(N/2)}\varphi(x_{i}, x_{j}) + h\,\psi(x_{i}, x_{j})
~.
\end{equation}
Note that $R_{|i-j|}^{(N/2)}$ is a dense circulant matrix, and all $N^2$ elements
differ from the standard Nystr\"om matrix \eqref{eq:standquad}.
Since $\varphi$ and $\psi$ do not usually have fast potential-theory based
algorithms to apply them, the Kress scheme is not FMM-compatible.

\section{Numerical experiments}
\label{sec:num}

We now describe numerical experiments in 1D, 2D and 3D applications
that illustrate the performance
of the quadratures described in sections \ref{sec:KR}-\ref{sec:kress}.
The experiments were carried out on a Macbook Pro with 2.4GHz Intel
Core 2 Duo and 4GB of RAM, and executed in a MATLAB environment.
Once the Nystr\"om matrix is filled, the linear system \eqref{eq:park4}
is solved via MATLAB's backslash ({\tt mldivide}) command.
In all examples below, the errors reported
are relative errors measured in the $L^{\infty}$-norm,
$||u_{\epsilon} - u||_{\infty}/||u||_{\infty}$, where $u$ is the reference solution
and $u_{\epsilon}$ is the numerical solution. For each experiment we compare the
performance of the different quadratures. Specifically, we compare the rules of
Kapur--Rokhlin of orders 2, 6, and 10; Alpert of orders 2, 6, 10, and 16;
modified Gaussian with $n=10$ points per panel; and (where convenient) Kress.
Our implementation of the modified Gaussian rule uses
$m=20$ auxiliary nodes for source and target on the same
panel, and $m'=24$ when on adjacent panels.

The quadrature nodes and weights used are provided in appendices and at the website \cite{2012_nystrom_website}.

\begin{figure}[htbp]
\centering
\includegraphics[width=0.65\textwidth]{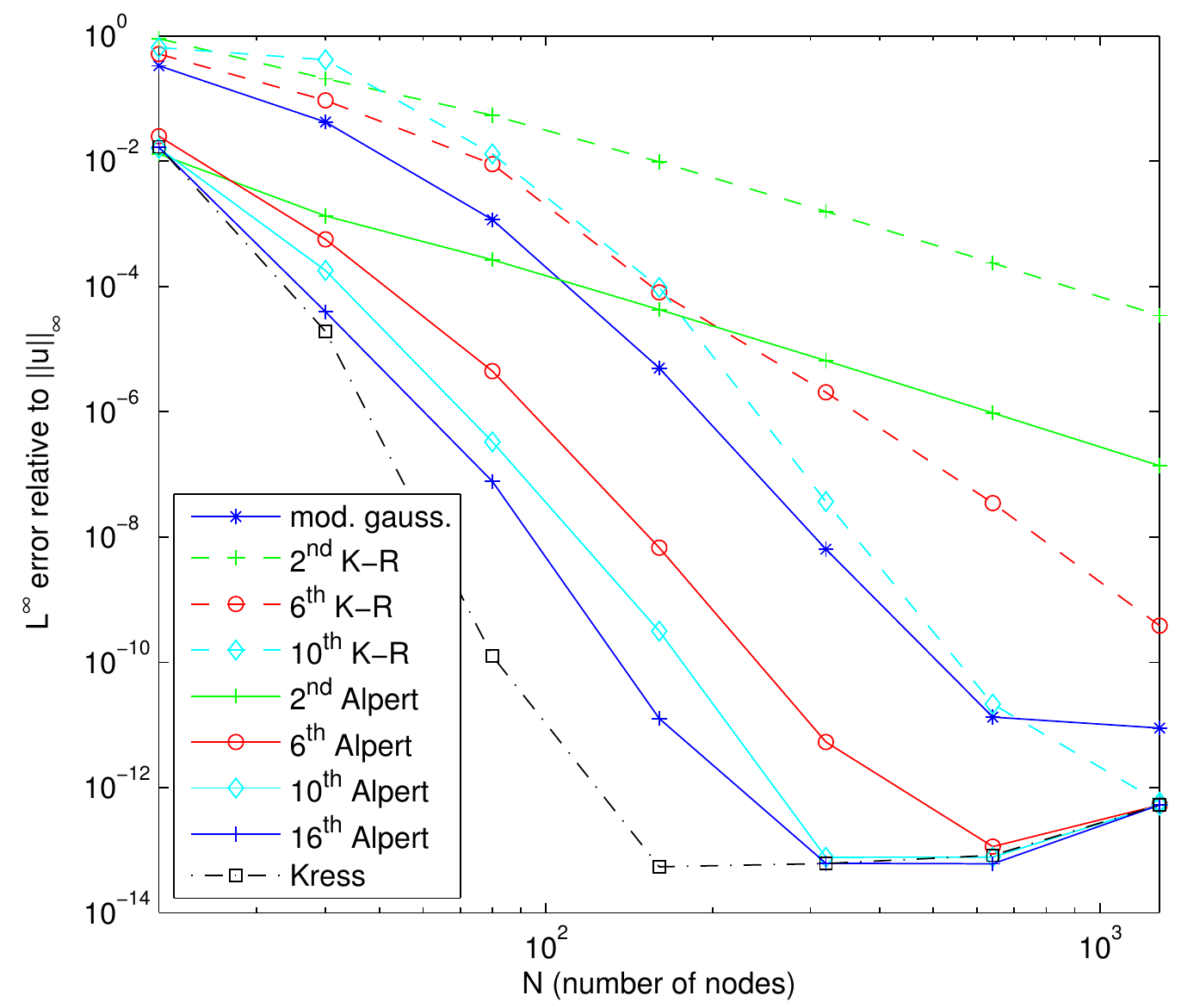}
\label{fig:1d}
\caption{Error results for solving the integral equation (\ref{eq:1d_eq}) in Section \ref{sec:1d_num}.}
\end{figure}

\subsection{A 1D integral equation example}
\label{sec:1d_num}
We solve the one-dimensional integral equation
\begin{equation}
\label{eq:1d_eq}
u(x) + \int_{0}^{2\pi}k(x,x')u(x')dx' = f(x),\qquad x \in [0,\,2\pi]
\end{equation}
associated with a simple kernel function having a periodic log
singularity at the diagonal,
\begin{equation}
k(x,x') = \frac{1}{2}\log\left|\sin\frac{x-x'}{2}\right| =
\frac{1}{4}\log \left(4\sin^2 \frac {t-s}{2}\right) - \frac{1}{2}\log 2~,
\end{equation}
thus the smooth
functions $\varphi(x,x') = 1/4$ and $\psi(x,x')=-(1/2)\log 2$ are constant.
This kernel is similar to that
arising from the Laplace single-layer operator on the unit circle.
Using \eqref{logfou} one may check that the above integral
operator has exact eigenvalues
$-\pi \log 2$ (simple) and $-\pi/(2n)$, $n=1,2,\ldots$ (each doubly-degenerate).
Thus the exact condition number of the problem \eqref{eq:1d_eq} is
$((\pi\log 2) -1)/(1-\pi/4) \approx 5.5$.
We choose the real-analytic periodic right-hand side
$f(x) = \sin(3x)\,e^{\cos(5x)}$.
The solution $u$ has $\|u\|_\infty \approx 6.1$.
We estimate errors by comparing to the
Kress solution at $N=2560$.
(In passing we note that the exact solution to \eqref{eq:1d_eq}
could be written analytically
as a Fourier series since the Fourier series of $f$ is
known in terms of modified Bessel functions.)
When a solution is computed on panel-based nodes, we
use interpolation back to the uniform trapezoid grid by evaluating
the Lagrange basis on the $n=10$ nodes on each panel.

In Figure \ref{fig:1d},
the errors in the $L^{\infty}$-norm divided by $\|u\|_\infty$
are presented for $N = 20, 40, 80, \dots, 1280$.
We see that the rules of order 2, 6, and 10 have the expected convergence
rates, but that Alpert has prefactors a factor $10^2$ to $10^5$ smaller
the Kapur--Rokhlin.
We also see that Kress is the most efficient at any desired accuracy,
followed by the three highest-order Alpert schemes.
These four schemes flatten out at 13 digits, but errors start
to grow again for larger $N$, believed due to the larger linear system.
Note that modified Gaussian performs roughly as well as 6th-order
Alpert with twice the number of points,
and that it flattens out at around 11 digits.

\begin{figure}[htbp]
\begin{minipage}{1.3in}(a) $0.5 \; \lambda$ diameter\\
\mbox{}\qquad($\omega=2.8$)
\end{minipage}
\quad \raisebox{-2.3in}%
{\includegraphics[width=0.6\textwidth]{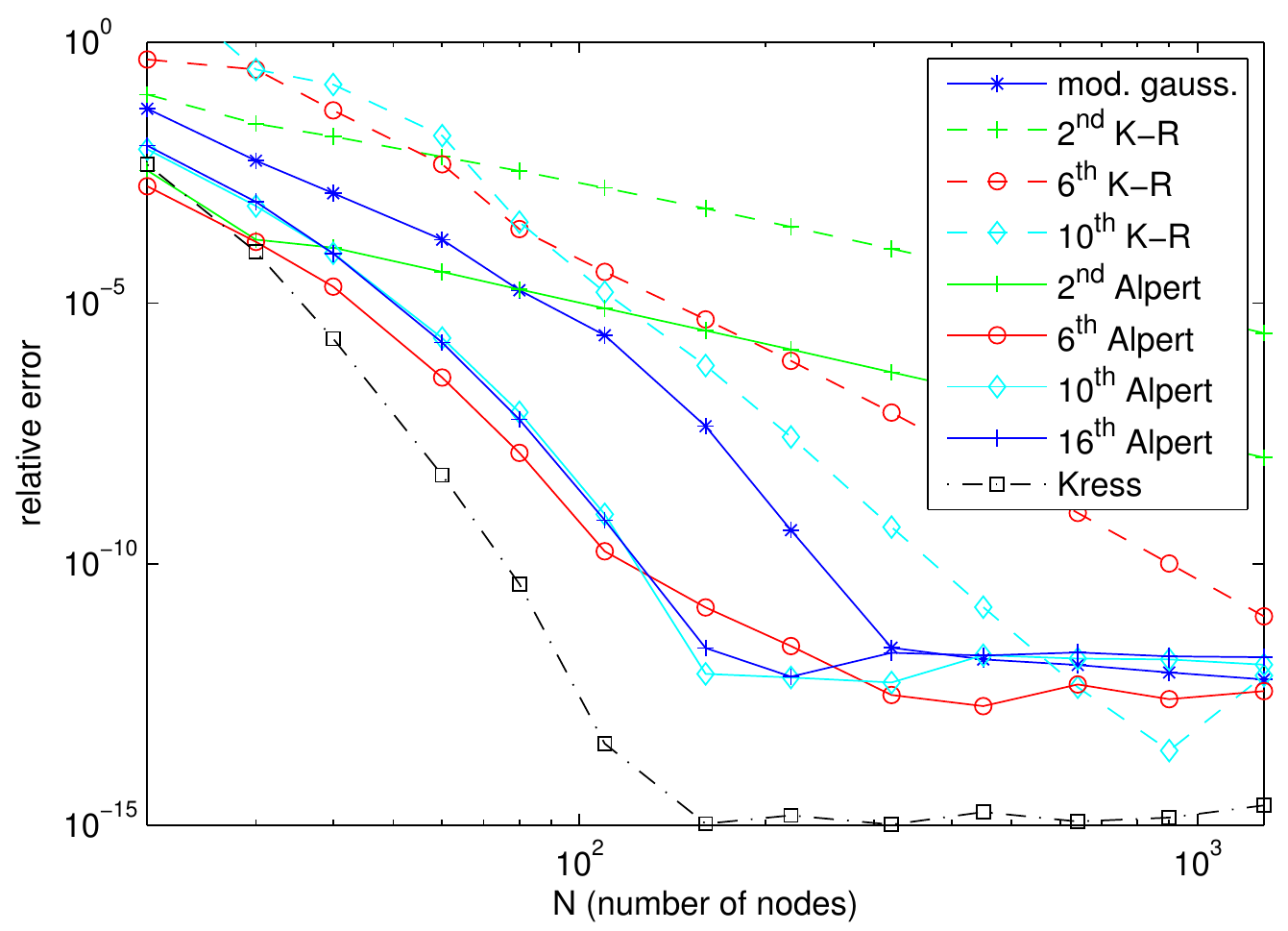}}
\\
\begin{minipage}{1.3in}(b) $5 \; \lambda$ diameter\\
\mbox{}\qquad($\omega=28$)
\end{minipage}
\quad\; \raisebox{-2.3in}%
{\includegraphics[width=0.6\textwidth]{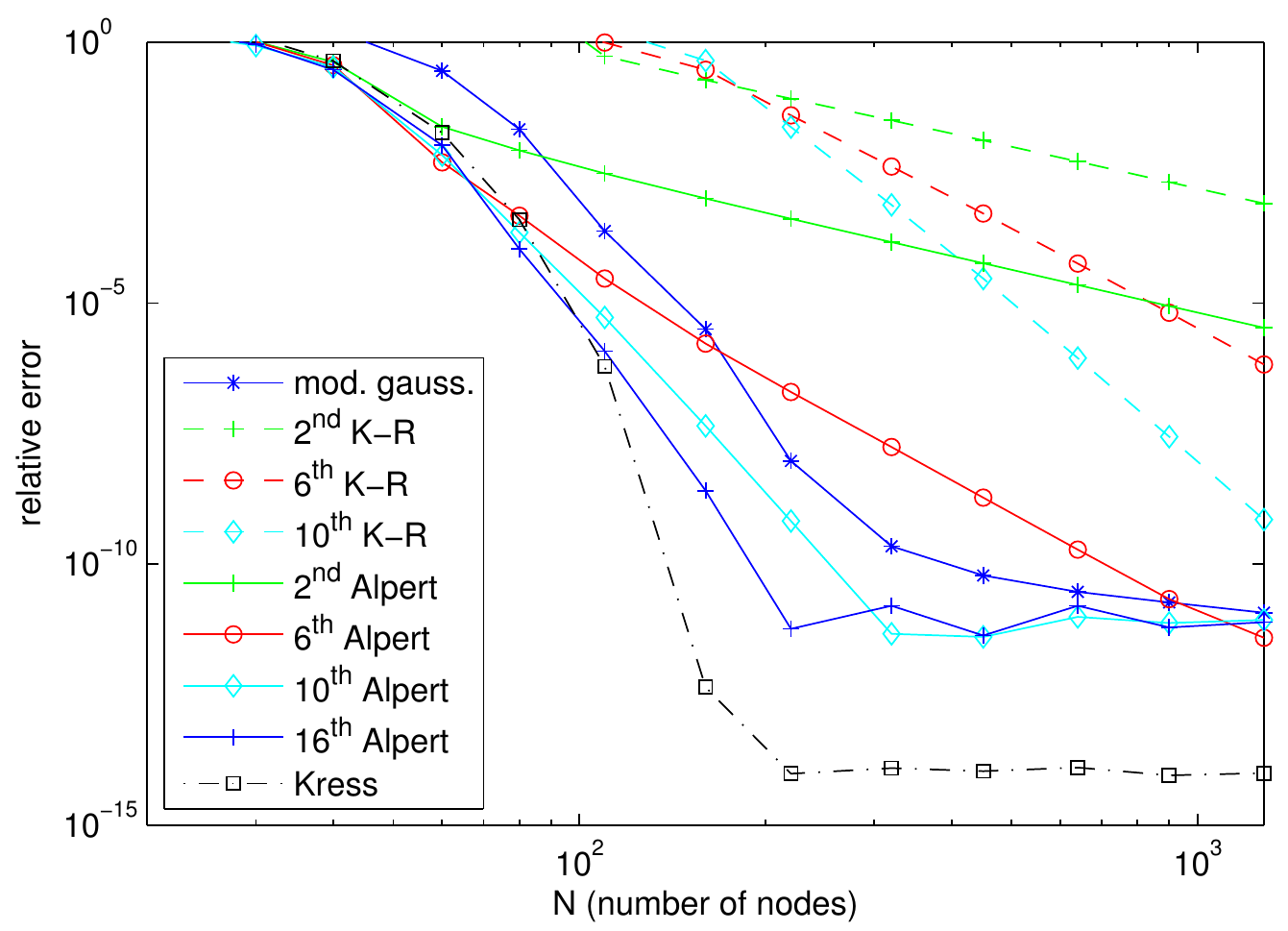}}
\\
\begin{minipage}{1.3in}(c) $50 \; \lambda$ diameter\\
\mbox{}\qquad($\omega=280$)
\end{minipage}
\quad \raisebox{-2.3in}%
{\includegraphics[width=0.6\textwidth]{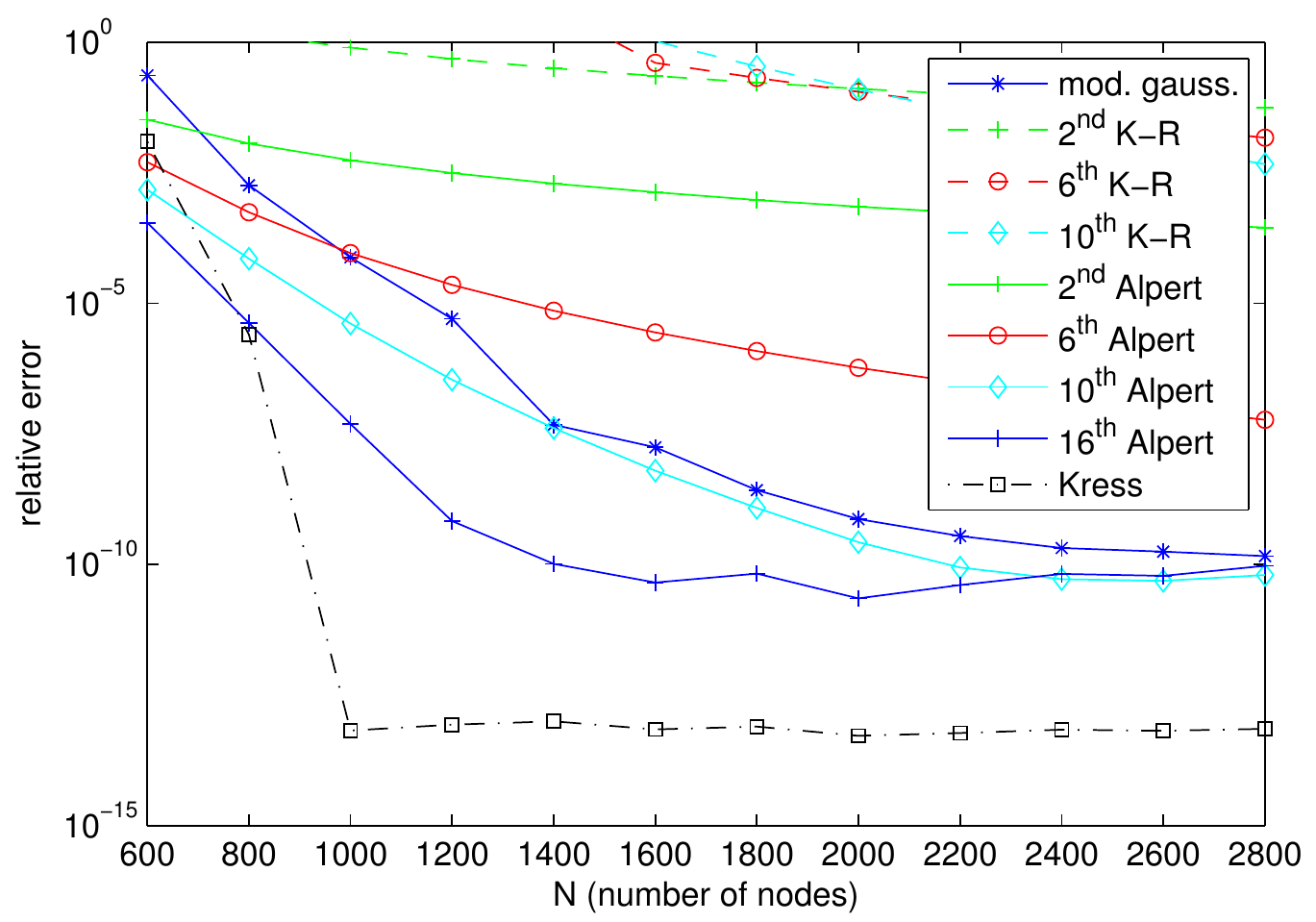}}
\caption[]{Error results for the exterior planar Helmholtz problem
\eqref{eq:helmholtz1} in Section \ref{sec:BIE_helm} solved on the
starfish domain of Figure \ref{fig:2Ddomains}.
\label{f:helm}
}
\end{figure}

\subsection{Combined field discretization of the Helmholtz equation in $\mathbb{R}^2$}
\label{sec:BIE_helm}
In this section, we solve the Dirichlet problem for the Helmholtz equation
exterior to a smooth domain $\Omega\subset\mathbb{R}^2$ with boundary $\Gamma$,
\begin{align}
\label{eq:helmholtz1}
-\Delta u - \omega^2 u &= 0, \qquad \text{in} \ \ E = \Omega^c, \\
u &= f, \qquad \text{on} \ \ \Gamma ,
\label{eq:helmholtz2}
\end{align}
where $\omega>0$ is the wavenumber. $u$ satisfies the Sommerfeld radiation condition
\begin{equation}
\lim_{r \rightarrow \infty} r^{1/2}
\left(\frac {\partial u}{\partial r} - i\omega u\right) = 0,
\end{equation}
where $r = |\vct{x}|$ and the limit holds uniformly in all directions.
A common approach \cite[Ch.~3]{coltonkress}
is to represent the solution to (\ref{eq:helmholtz1}) via both the single and double layer acoustic potentials,
\begin{align}
\label{eq:sing_doub_represent}
u(\vct{x}) &= \int_\Gamma k(\vct{x}, \vct{x}')\,\sigma(\vct{x}')\,dl(\vct{x}') \nonumber \\
       &= \int_\Gamma
       \left(\frac {\partial \phi(\vct{x}, \vct{x}')}{\partial \vct{n}(\vct{x}')} - i\omega\,\phi(\vct{x},\vct{x}')\right)\,
       \sigma(\vct{x}')\,dl(\vct{x}'),  \qquad \vct{x}\in E,
\end{align}
where $\phi(\vct{x}, \vct{x}') = \frac {i}{4}H_{0}^{(1)}(\omega|\vct{x} - \vct{x}'|)$ and $H_{0}^{(1)}$ is the Hankel function of the first kind of order zero; $\vct{n}$ is the normal vector pointing outward to $\Gamma$, and
$dl$ the arclength measure on $\Gamma$. The motivation for the
combined representation  (\ref{eq:sing_doub_represent}) is to obtain the unique solvability to problem (\ref{eq:helmholtz1}-\ref{eq:helmholtz2}) for all $\omega>0$.
The corresponding boundary integral equation we need to solve is
\begin{equation}
\label{eq:BIE_helmholtz}
\frac{1}{2}\sigma(\vct{x})+\int_\Gamma k(\vct{x}, \vct{x}')\sigma(\vct{x}')\,dl(\vct{x}')  = f(\vct{x}), \qquad \vct{x}\in \Gamma,
\end{equation}
where $k(\vct{x}, \vct{x}') = \frac {\partial \phi(\vct{x}, \vct{x}')}{\partial \vct{n}(\vct{x}')} - i\omega\,\phi(\vct{x},\vct{x}')$.

To convert \eqref{eq:BIE_helmholtz} into an integral equation on the
real line, we need to parametrize $\Gamma$ by a
vector-valued smooth function $\vct{\tau} : [0,T] \rightarrow \mathbb{R}^2$.
By changing variable, \eqref{eq:BIE_helmholtz} becomes
\begin{equation}
\label{eq:1Dparam}
\frac{1}{2}\sigma(\vct{\tau}(t)) + \int_{0}^{T} k(\vct{\tau}(t),
\vct{\tau}(s))\,\sigma(\vct{\tau}(s))\,|d\vct{\tau}/ds|\,ds =
f(\vct{\tau}(t)),\qquad t \in [0,T].
\end{equation}
To keep our formula uncluttered, we rewrite the kernel as
\begin{equation}
m(t, s) =  k(\vct{\tau}(t), \vct{\tau}(s))\,|d\vct{\tau}/ds|,
\label{eq:simp_ker}
\end{equation}
as well as the functions
\begin{equation*}
\sigma(t) = \sigma(\vct{\tau}(t)) \hspace{1em} \textrm{and} \hspace{1em} f(t) = f(\vct{\tau}(t)).
\end{equation*}
Thus we may write the integral equation in standard form,
\begin{equation}
\sigma(t) + 2\int_{0}^{T}m(t,s)\,\sigma(s)\,ds  = 2f(s), \qquad t \in [0,T],
\label{eq:eq1D}
\end{equation}
and apply the techniques of this paper to it.

\begin{remark}
All the quadrature schemes apart from that of Kress are now easy to
implement by evaluation of $m(t,s)$.
However, to implement the Kress scheme,
separation of the parametrized single- and double-layer Helmholtz kernels
into analytic functions $\varphi$ and $\psi$ is necessary, and not trivial.
We refer the reader to \cite[Sec.~2]{kress91} or \cite[Ch.~3]{coltonkress}.
\end{remark}

We assess the accuracy of each quadrature rule for the smooth
domain shown in Figure \ref{fig:2Ddomains}(c),
varying $N$ and the wavenumber $\omega$.
Specifically, we varied wavenumbers such that there are $0.5$, $5$ and $50$
wavelengths across the domain's diameter.
The right-hand side is generated by a sum of five point sources
inside $\Omega$, with random strengths; thus the exact exterior
solution is known.
Errors are taken to be the maximum relative error
at the set of measurement points shown in
Figure \ref{fig:2Ddomains}(c).
Notice that sources and measurement points are both far from $\Gamma$,
thus no challenges due to close evaluation arise here.

The results are shown in Figure \ref{f:helm}.
At the lowest frequency, results are similar to Figure \ref{fig:1d},
except that errors bottom out at slightly higher accuracies,
probably because of the smoothing effect of evaluation of $u$
at distant measurement points.
In terms of the error level each scheme saturates at,
Kress has a 2-3 digit advantage over the others at all frequencies.
At the highest frequency, Figure \ref{f:helm}(c),
there are about 165 wavelengths around the perimeter $\Gamma$,
hence the $N=1000$ at which Kress is fully converged
corresponds to around 6 points per wavelength.
At 10 points per wavelength (apparently a standard choice in the engineering
community, and roughly halfway along the $N$ axis in our plot),
16th-order Alpert has converged at 10 digits,
while modified Gaussian and 10th-order Alpert are similar at 8 digits.
The other schemes are not competitive.

\begin{table}
{\scriptsize
\begin{tabular}{l|lllllllll}
scheme &mod. Gauss &2nd K-R &6th K-R &10th K-R &2nd Alpert &6th Alpert &
10th Alpert&16th Alpert& Kress\\
\hline
cond \# & 3.95 & 3.52 & 3.68 & 169 & 3.52 & 3.52 & 3.52 & 3.52 & 3.52\\
\# iters & 14&14&22&206&14&14&14&14&14\\
\hline
\end{tabular}
}
\vspace{2ex}
\caption{Condition numbers of the Nystr\"om system matrix $(\frac{1}{2}\mtx{I} + \mtx{A})$,
and numbers of GMRES iterations to reach residual error $10^{-12}$,
for all the quadrature schemes.
\label{t:GMRES}
}
\end{table}

\begin{figure}
(a)\raisebox{-1in}{\includegraphics[width=6in]{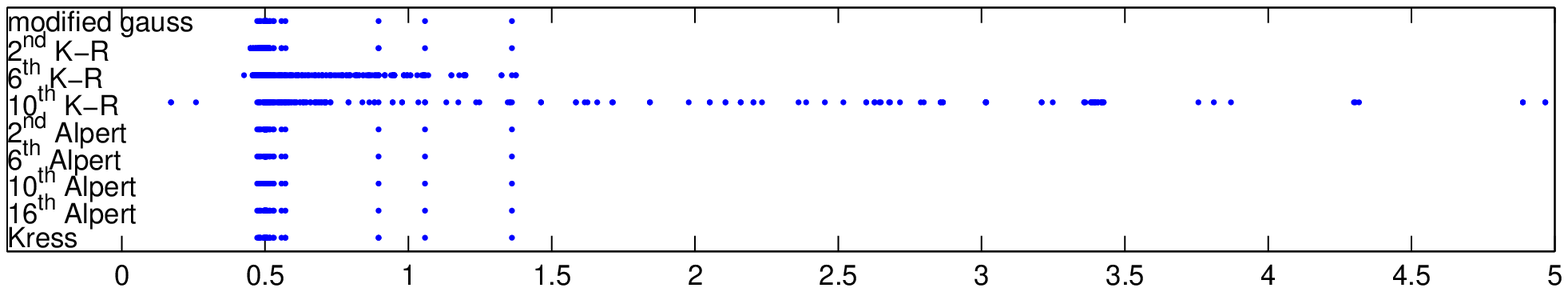}}
\\
(b)\raisebox{-1.8in}{\includegraphics[height=2in]{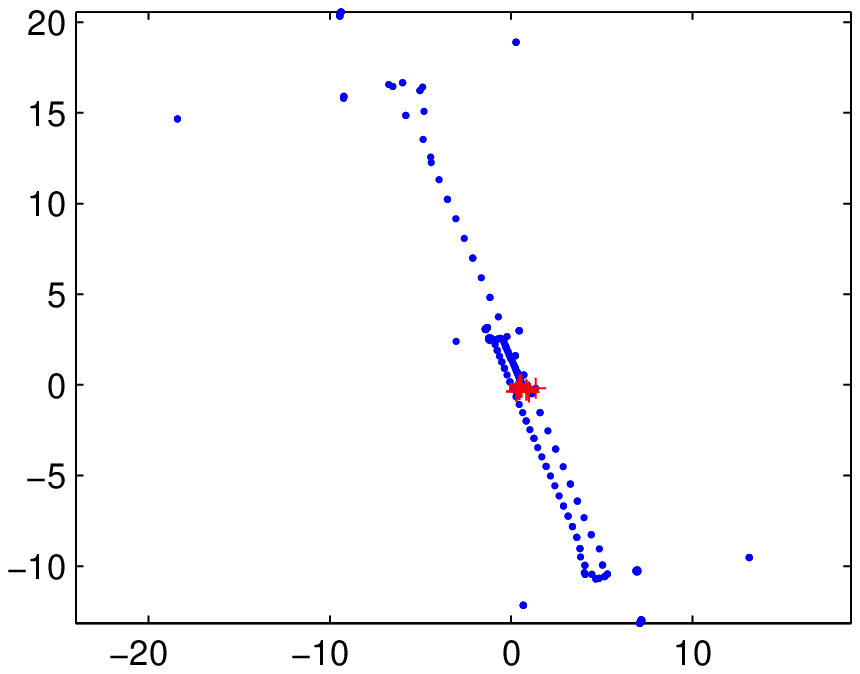}}
\qquad
(c)\raisebox{-1.8in}{\includegraphics[height=2in]{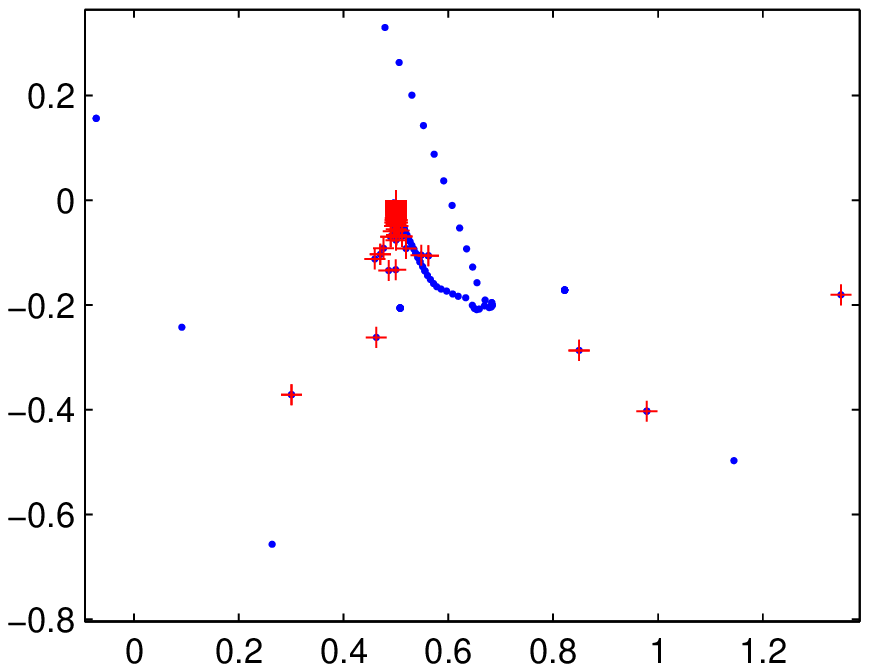}}
\caption{
(a) Magnitude of eigenvalues of the matrix $(\tfrac{1}{2}\mtx{I} + \mtx{A})$
associated with the Nystr\"om discretization of the Helmholtz BVP
(\ref{eq:BIE_helmholtz}). The system size is $N=640$ and the wave-number
$\omega$ corresponds to a contour of size $0.5$ wave-lengths.
(b)  Eigenvalues in the complex plane associated with 10th-order
Kapur--Rokhlin (dots) and Kress (crosses) quadratures.
(c) Same plot as (b), but zoomed in to the origin.
\label{f:eig}
}
\end{figure}

\subsection{Effect of quadrature scheme on iterative solution efficiency}
\label{sec:GMRES}

When the number of unknowns $N$ becomes large, iterative solution
becomes an important tool.
One standard iterative method for nonsymmetric systems is GMRES \cite{saad}.
In Table~\ref{t:GMRES} we compare the numbers of GMRES
iterations needed to solve the linear system \eqref{eq:park2}
arising from the low-frequency Helmholtz BVP
in section~\ref{sec:BIE_helm} (0.5 wavelengths across),
when the matrix was constructed via the various quadrature schemes.
We also give the condition number of the system.
We see that most schemes result in 14 iterations
and a condition number of around 3.5; this reflects
the fact that the underlying integral equation is Fredholm 2nd-kind
and well-conditioned.
However, 6th-order and particularly 10th-order Kapur--Rokhlin
require many more iterations (a factor 15 more in the latter case),
and have correspondingly high condition numbers.
In a practical setting this would mean a much longer solution time.

In order to understand why, in Figure \ref{f:eig}
we study the spectra of the system matrix; since the operator is of the form
Id$/2$ + compact, eigenvalues should cluster near 1/2.
Indeed this is the case for all the schemes.
However, 6th- and 10th-order Kapur--Rokhlin also contain additional
eigenvalues that do not appear to relate to those of the operator.
In the 10th-order case, these create a wide
``spray'' of many eigenvalues,
which we also plot in the complex plane in (b) and (c).
A large number (around 200) of eigenvalues fall at much larger distances
than the true spectrum,
and on both sides of the origin; we believe they cause the slow
GMRES convergence.
The corresponding eigenvectors are oscillatory, typically
alternating in sign from node to node.
We believe this pollution of the spectrum arises
from the large alternating weights $\gamma_l$ in these schemes.
Note that one spurious eigenvalue falls near the origin;
this mechanism could induce an arbitrarily large condition number
even though the integral equation condition number is small.
Although we do not test the very large $N$ values where iterative methods
become essential, we expect that our conclusions apply also to larger $N$.

\begin{figure}[ht!]
\begin{center}
\begin{minipage}{0.31\linewidth} \begin{flushleft}
\hspace{10mm}
\includegraphics[height=.75\linewidth]{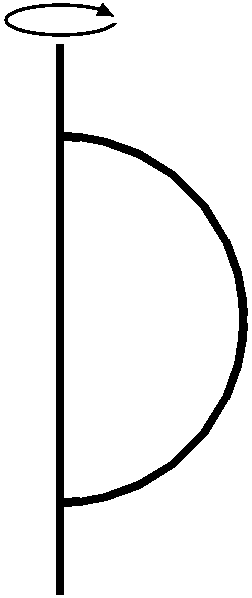}
\end{flushleft} \end{minipage}
\begin{minipage}{0.31\linewidth} \begin{center}
\includegraphics[height =.7\linewidth]{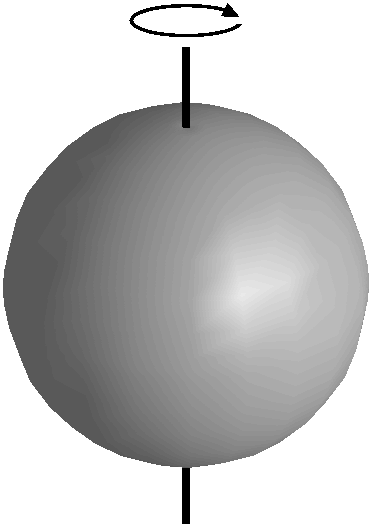}
\end{center} \end{minipage}
\\
\begin{minipage}{1\linewidth}\begin{center} (a) \end{center} \end{minipage}
\\ \vspace{4mm}
\begin{minipage}{0.31\linewidth} \begin{flushleft}
\hspace{10mm}
\includegraphics[height=.75\linewidth]{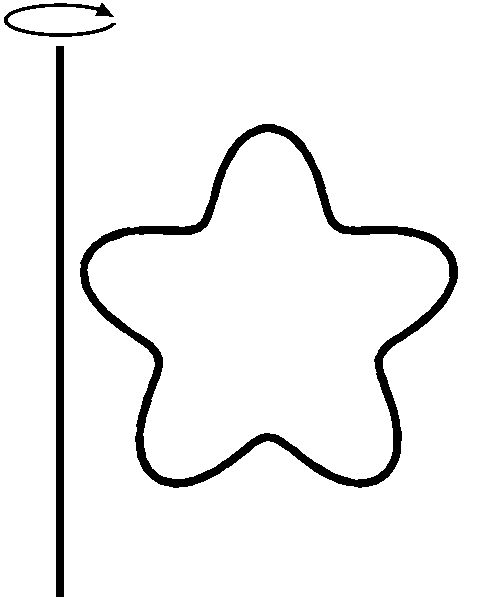}
\end{flushleft} \end{minipage}
\begin{minipage}{0.31\linewidth} \begin{center}
\includegraphics[height =.55\linewidth]{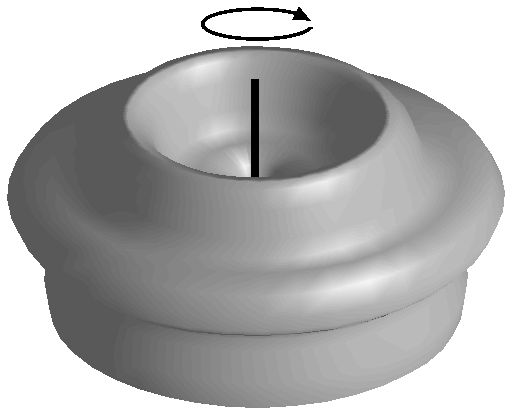}
\end{center} \end{minipage}
\\
\begin{minipage}{1\linewidth}\begin{center} (b) \end{center} \end{minipage}
\\ \vspace{4mm}
\caption{Domains used in numerical examples in Section \ref{sec:3dBIE}. All items are rotated about the vertical axis.  (a) A sphere.   (b) A starfish torus.  }
\label{fig:3Ddomains}
\end{center}
\end{figure}

\begin{figure} 
(a)\raisebox{-3.4in}{\includegraphics[width=0.65\textwidth]{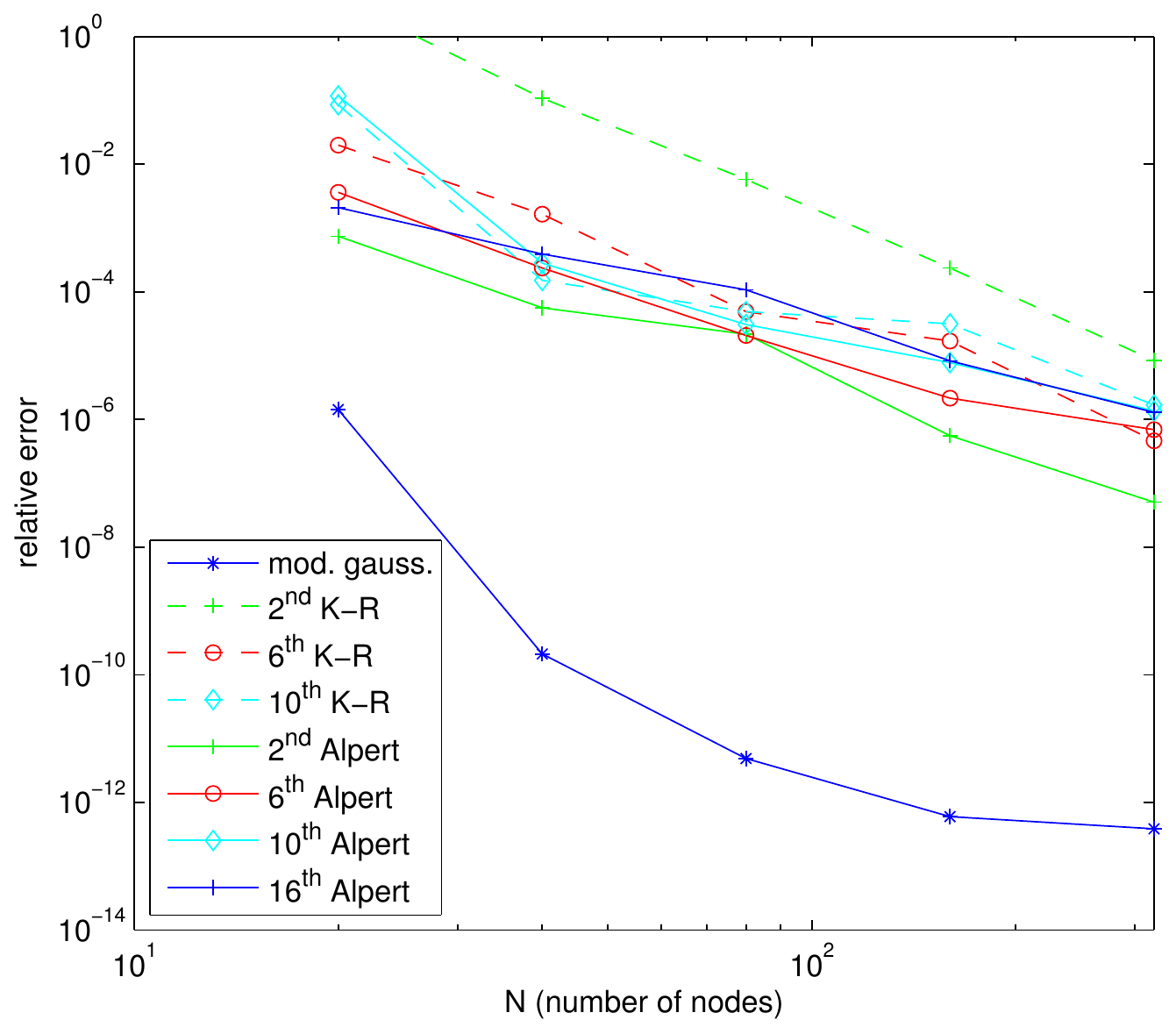}}
\\
(b)\raisebox{-3.4in}{\includegraphics[width=0.65\textwidth]{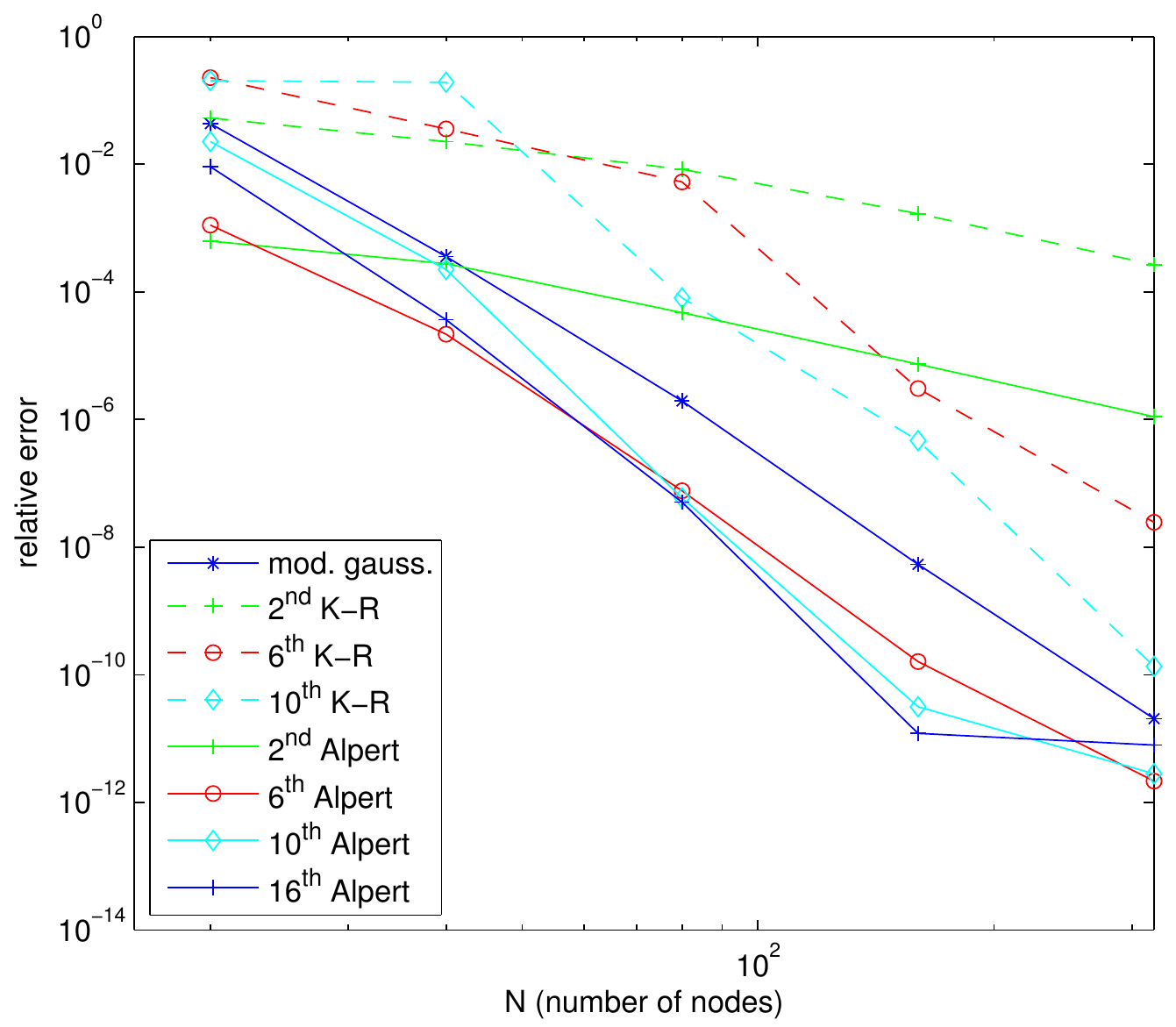}}
\caption{Error results for the 3D
interior Dirichlet Laplace problem from
section \ref{sec:3dBIE} solved on the axisymmetric domains (a) and (b)
respectively
shown in Figure \ref{fig:3Ddomains}.
\label{f:3d}}
\end{figure}

\subsection{The Laplace BVP on axisymmetric surfaces in $\mathbb{R}^{3}$}
\label{sec:3dBIE}

In this section, we compare quadratures rules applied on kernels associated with BIEs on rotationally symmetric surfaces in $\mathbb{R}^{3}$. Specifically, we considered integral equations of the form
\begin{equation}
\label{eq:axisymmetric_3d}
\sigma(\vct{x}) + \int_\Gamma k(\vct{x}, \vct{x}')\,\sigma(\vct{x}')\,dA(\vct{x}') = f(\vct{x}),  \quad x\in \Gamma,
\end{equation}
under the assumptions that $\Gamma$ is a surface in $\mathbb{R}^{3}$ obtained by rotating a curve $\gamma$ about an axis and the kernel function $k$ is invariant under rotation about the symmetry axis.  Figure \ref{fig:3Ddomains} depicts domains used in numerical examples: the generating curves $\gamma$ are shown in the left figures and the axisymmetric surfaces $\Gamma$ are shown in the right ones. The BIE (\ref{eq:axisymmetric_3d}) on rotationally symmetric surfaces can via a Fourier transform be recast as a sequence of equations defined on the generating curve in cylindrical coordinates, i.e.
\begin{equation}
\label{eq:axisymmetric_2d}
\sigma_n(r, z) + \sqrt{2\pi}\,\int_\Gamma k_n(r, z, r', z')\,\sigma_n(r', z')\,r'\,dl(r', z') = f_n(r, z), \quad (r, z)\in \gamma, \quad n\in\mathbb{Z},
\end{equation}
where $\sigma_n$, $f_n$, and $k_n$ denote the Fourier coefficients of $\sigma$, $f$, and $k$, respectively. Details on how to truncate the Fourier series and construct the coefficient matrices for Laplace problem and Helmholtz problem can be found in \cite{young12}. In the following experiments, we consider the BIE (\ref{eq:axisymmetric_3d}) which arises from the interior Dirichlet Laplace problem, in which case
\begin{equation}
k(\vct{x}, \vct{x}') = \frac {\vct{n}(\vct{x}')\cdot(\vct{x}-\vct{x'})}{4\pi|\vct{x}-\vct{x'}|^3}.
\end{equation}
As we recast the BIE defined on $\Gamma$ to a sequence of equations defined on the
generating curve $\gamma$, it is easy to see that the kernel function $k_{n}$ has a
logarithmic singularity as $(r',z') \rightarrow (r,z)$.
In this experiment, $101$ Fourier modes were used.
We tested all the quadrature schemes apart from that of Kress, since
we do not know of an analytic split of the axisymmetric kernel $k_n$ into
smooth parts $\varphi$ and $\psi$.


Equation (\ref{eq:axisymmetric_3d}) was solved for Dirichlet data $f$ corresponding
to an exact solution $u$ generated by point charges placed outside the domain.
The errors reported reflect the maximum of the point-wise errors (compared to
the known analytic solution) sampled at a set of target points inside the domain.

The results are presented in Figure \ref{f:3d}.
The most apparent feature in (a) is that, because the curve $\gamma$ is
open, the schemes based on the periodic trapezoid rule fail to give
high-order convergence; rather, it appears to be approximately 3rd-order.
Panel-based schemes are able to handle open intervals as easily as periodic
ones, thus modified Gaussian performs well: it reaches 12-digit accuracy with
only around $N=100$ points.
(We remark that the problems associated with an open curve are in this case
artificial and can be overcome by a better problem formulation. We deliberately
chose to use a simplistic formulation to simulate an ``open curve'' problem.)
In (b), all functions are again periodic since $\gamma$ is closed;
modified Gaussian performs similarly to the three highest-order
Alpert schemes with around 1.5 to 2 times the number of points.

\section{Concluding remarks}
\label{sec:conc}

To conclude, we make some informal remarks on the relative advantages
and disadvantages of the different quadrature rules that we have
discussed. Our remarks are informed primarily by the numerical
experiments in section \ref{sec:num}.

Comparing the three schemes based upon nodes equi-spaced in parameter
(Kapur--Rokhlin, Alpert, and Kress),
we see that Kress always excels due to its superalgebraic convergence,
converging fully at around 6 points per wavelength at high frequency,
and its small saturation error of $10^{-13}$ to $10^{-15}$.
However, the analytic split required for Kress is not always available
or convenient, and Kress is not amendable to standard FMM-style
fast matrix algebra. Kapur--Rokhlin and Alpert show their expected
algebraic convergence rates at orders 2, 6, and 10,
and both seem to saturate at around $10^{-12}$.
However, Alpert outperforms
Kapur--Rokhlin since its
prefactor is much lower, resulting in 2-8 extra digits of accuracy at
the same $N$.
Another way to compare these two is that Kapur--Rokhlin requires
around 6 to 10 times the number of unknowns as Alpert to
reach comparable accuracy.
The difference in a high-frequency problem is striking, as in
Figure \ref{f:helm}.
The performance of 16th-order Alpert is usually less than 1 digit
better than 10th-order Alpert,
apart from at high frequency when it can be up to 3 digits better.

Turning to the panel-based modified Gaussian scheme,
we see that in the low-frequency settings it behaves
like 10th-order Alpert but requires around 1.5 to 2 times the $N$ to
reach similar accuracy.
This may be related to the fact that Gauss-Legendre panels would need a factor
$\pi/2$ higher $N$ than the trapezoid rule
to achieve the same largest spacing between nodes;
this is the price to pay for a panel-based scheme.
However, for
medium and high frequencies, $10$th-order Alpert has
little distinguishable advantage over modified Gaussian.
Both are able to reach around 10 digit accuracy at 15 points per wavelength.
Modified Gaussian errors seem to saturate at around $10^{-11}$ to $10^{-12}$.
It therefore provides a good all-round choice, expecially when
adaptivity is anticipated, or global parametrizations are not readily constructed.
One disadvantage relative to the other schemes
is that the auxiliary nodes require kernel
evaluations that are very close to the singularity ($10^{-7}$ or less;
for Alpert the minimum is only around $10^{-3}$).

We have not tested situations in which adaptive quadrature
becomes essential; in such cases modified Gaussian would excel.
However, a hint of the convenience of modified Gaussian is
given by its effortless handling of an open curve in Figure \ref{f:3d}(a)
where the other (periodic) schemes become low-order
(corner-style reparametrizations would be needed to fix this
\cite[Sec.~3.5]{coltonkress}).

In addition, we have showed that, in an iterative solution setting,
higher-order Kapur--Rokhlin
can lead to much slower GMRES convergence than any of the other schemes.
We believe this is because it introduces many large eigenvalues into the
spectrum, unrelated to those of the underlying operator.
Thus 10th-order Kapur--Rokhlin should be used with caution.
With that said, Kapur--Rokhlin is without doubt
the simplest to implement of the four schemes, since no interpolation
or new kernel evaluations are needed.

We have chosen to not report computational times in this note since our
MATLAB implementations are far from optimized for speed. However, it should
be mentioned that both the Alpert method and the method based on modified
Gaussian quadrature require a substantial number (between $20N$ and $30N$ in
the examples reported) of additional kernel evaluations.

%

For simplicity,
in this note we limited our attention to the case of smooth contours,
but both the Alpert and the modified Gaussian rule can with certain modifications be
applied to contours with corners, see, e.g.,
\cite{2011_ojala_diss,Bremer:a,Bremer:b,2011_bremer,Helsing:00a,helsingtut}.
We plan to include
the other recent schemes shown in
Table~\ref{t:schemes}, and curves with corners, in future comparisons.

\section*{Acknowledgments}
We acknowledge the use of a small quadrature code by Andras Pataki.
We are grateful for useful discussions with Andreas Kl\"ockner about operator eigenvalues.
The work of SH and PGM is supported by NSF grants DMS-0748488 and CDI-0941476.
The work of AHB is supported by NSF grant DMS-1216656.

\bibliographystyle{siam}
\bibliography{refs}

\begin{thebibliography}{10}

\bibitem{alpert:99a}
{\sc B.~K. Alpert}, {\em Hybrid gauss-trapezoidal quadrature rules}, SIAM J.
  Sci. Comput., 20 (1999), pp.~1551--1584.

\bibitem{atkinson1997}
{\sc K.~E. Atkinson}, {\em The numerical solution of integral equations of the
  second kind}, Cambridge University Press, Cambridge, 1997.

\bibitem{1986_barnes_hut}
{\sc J.~Barnes and P.~Hut}, {\em A hierarchical $o(n \log n)$ force-calculation
  algorithm}, Nature, 324 (1986).

\bibitem{2011_bremer}
{\sc J.~Bremer}, {\em A fast direct solver for the integral equations of
  scattering theory on planar curves with corners}, Journal of Computational
  Physics,  (2011), pp.~--.

\bibitem{Bremer:a}
{\sc J.~Bremer and V.~Rokhlin}, {\em Efficient discretization of laplace
  boundary integral equations on polygonal domains}, J. Comput. Phys., 229
  (2010), pp.~2507--2525.

\bibitem{Bremer:b}
{\sc J.~Bremer, V.~Rokhlin, and I.~Sammis}, {\em Universal quadratures for
  boundary integral equations on two-dimensional domains with corners}, Journal
  of Computational Physics, 229 (2010), pp.~8259 -- 8280.

\bibitem{coltonkress}
{\sc D.~Colton and R.~Kress}, {\em Inverse acoustic and electromagnetic
  scattering theory}, vol.~93 of Applied Mathematical Sciences,
  Springer-Verlag, Berlin, second~ed., 1998.

\bibitem{2001_krasny_treecode}
{\sc Z.-H. Duan and R.~Krasny}, {\em An adaptive treecode for computing
  nonbonded potential energy in classical molecular systems}, Journal of
  Computational Chemistry, 22 (2001), pp.~184--195.

\bibitem{rokhlin1987}
{\sc L.~Greengard and V.~Rokhlin}, {\em A fast algorithm for particle
  simulations}, J. Comput. Phys., 73 (1987), pp.~325--348.

\bibitem{2012_nystrom_website}
{\sc S.~Hao, A.~Barnett, and P.~Martinsson}, {\em Nystr\"om quadratures for
  {BIE}s with weakly singular kernels on 1{D} domains}, 2012.
\newblock http://amath.colorado.edu/faculty/martinss/Nystrom/.

\bibitem{helsingmixed}
{\sc J.~Helsing}, {\em Integral equation methods for elliptic problems with
  boundary conditions of mixed type}, J. Comput. Phys., 228 (2009),
  pp.~8892--8907.

\bibitem{helsingtut}
\leavevmode\vrule height 2pt depth -1.6pt width 23pt, {\em Solving integral
  equations on piecewise smooth boundaries using the {RCIP} method: a
  tutorial}, 2012.
\newblock preprint, 34 pages, {\tt arXiv:1207.6737v3}.

\bibitem{Helsing:00a}
{\sc J.~Helsing and R.~Ojala}, {\em Corner singularities for elliptic problems:
  Integral equations, graded meshes, quadrature, and compressed inverse
  preconditioning}, J. Comput. Phys., 227 (2008), pp.~8820--8840.

\bibitem{Kapur:97a}
{\sc S.~Kapur and V.~Rokhlin}, {\em High-order corrected trapezoidal quadrature
  rules for singular functions}, SIAM J. Numer. Anal., 34 (1997),
  pp.~1331--1356.

\bibitem{qbx}
{\sc A.~Kl\"ockner, A.~H. Barnett, L.~Greengard, and M.~O'Neil}, {\em
  Quadrature by expansion: a new method for the evaluation of layer
  potentials}, 2012.
\newblock submitted.

\bibitem{Kolm:01a}
{\sc P.~Kolm and V.~Rokhlin}, {\em Numerical quadratures for singular and
  hypersingular integrals}, Comput. Math. Appl., 41 (2001), pp.~327--352.

\bibitem{kress91}
{\sc R.~Kress}, {\em Boundary integral equations in time-harmonic acoustic
  scattering}, Mathl. Comput. Modelling, 15 (1991), pp.~229--243.

\bibitem{LIE}
{\sc R.~Kress}, {\em Linear Integral Equations}, vol.~82 of Applied
  Mathematical Sciences, Springer, second~ed., 1999.

\bibitem{Martinsson:04a}
{\sc P.~Martinsson and V.~Rokhlin}, {\em A fast direct solver for boundary
  integral equations in two dimensions}, J. Comput. Phys., 205 (2004),
  pp.~1--23.

\bibitem{nystrom}
{\sc E.~Nystr{\"o}m}, {\em {\"U}ber die praktische {Aufl\"osung} von
  {Integralgleichungen} mit {Andwendungen} aug {Randwertaufgaben}}, Acta Math.,
  54 (1930), pp.~185--204.

\bibitem{2011_ojala_diss}
{\sc R.~Ojala}, {\em Towards an All-Embracing Elliptic Solver in 2{D}}, PhD
  thesis, Department of Mathematics, Lund University, Sweden, 2011.

\bibitem{saad}
{\sc Y.~Saad}, {\em Iterative Methods for Sparse Linear Systems}, Society for
  Industrial and Applied Mathematics, 2nd ed.~ed., 2003.

\bibitem{tref}
{\sc L.~N. Trefethen}, {\em Spectral methods in {MATLAB}}, vol.~10 of Software,
  Environments, and Tools, Society for Industrial and Applied Mathematics
  (SIAM), Philadelphia, PA, 2000.

\bibitem{ATAP}
{\sc L.~N. Trefethen}, {\em Approximation Theory and Approximation Practice},
  SIAM, 2012.
\newblock {\tt http://www.maths.ox.ac.uk/chebfun/ATAP}.

\bibitem{kifmm}
{\sc L.~Ying, G.~Biros, and D.~Zorin}, {\em A kernel-independent adaptive fast
  multipole method in two and three dimensions}, J. Comput. Phys., 196 (2004),
  pp.~591--626.

\bibitem{young12}
{\sc P.~M. Young, S.~Hao, and P.~G. Martinsson}, {\em A high-order {Nystr\"om}
  discretization scheme for boundary integral equations defined on rotationally
  symmetric surfaces}, J. Comput. Phys., 231 (2012), pp.~4142--4159.

\end{thebibliography}

\appendix

\section{Tables of Kapur--Rokhlin Quadrature weights}
\label{a:KRrules}

{\footnotesize 

\begin{center}

\begin{minipage}{0.48\linewidth} \begin{center}
\begin{tabular}{ | c | c | }
\hline
  \multicolumn{2}{|c|}{$2^{nd}$-order Kapur--Rokhlin correction weights}   \\
  \multicolumn{2}{|c|}{integrals of the form $\int_{0}^{1}f(x)+g(x)\log(x)\, dx$}\\
\hline
    INDEX $l$ & WEIGHTS $\gamma_l + \gamma_{-l}$ \\
\hline
   1 &  1.825748064736159e+00\\
2 & -1.325748064736159e+00 \\
\hline
\end{tabular}
\end{center} \end{minipage}\\

\vspace{3em}
\begin{minipage}{0.48\linewidth} \begin{center}
\begin{tabular}{ | c | c | }
\hline
  \multicolumn{2}{|c|}{$6^{th}$-order Kapur--Rokhlin correction weights}   \\
  \multicolumn{2}{|c|}{integrals of the form $\int_{0}^{1}f(x)+g(x)\log(x)\, dx$}
\\
\hline
    INDEX $l$ & WEIGHTS $\gamma_l + \gamma_{-l}$\\
\hline
     1&4.967362978287758e+00\\ 2&-1.620501504859126e+01\\ 3&2.585153761832639e+01 \\
        4&-2.222599466791883e+01\\ 5&9.930104998037539e+00\\ 6&-1.817995878141594e+00\\
\hline
\end{tabular}
\end{center} \end{minipage}\\

\vspace{3em}
\begin{minipage}{0.48\linewidth} \begin{center}
\begin{tabular}{ | c | c | }
\hline
  \multicolumn{2}{|c|}{$10^{th}$-order Kapur--Rokhlin correction weights}   \\
  \multicolumn{2}{|c|}{integrals of the form $\int_{0}^{1}f(x)+g(x)\log(x)\, dx$}
\\
\hline
    INDEX $l$ & WEIGHTS $\gamma_l + \gamma_{-l}$\\
\hline
1&7.832432020568779e+00\\ 2&-4.565161670374749e+01 \\3&1.452168846354677e+02 \\
        4&-2.901348302886379e+02 \\ 5& 3.870862162579900e+02 \\6&-3.523821383570681e+02 \\
        7&2.172421547519342e+02\\ 8& -8.707796087382991e+01 \\ 9&2.053584266072635e+01 \\
       10&-2.166984103403823e+00\\
\hline
\end{tabular}
\end{center} \end{minipage}\\
\end{center}

} 

\section{Tables of Alpert Quadrature rules}
\label{a:Alprules}

\begin{center}
\begin{minipage}{0.48\linewidth} \begin{center}
\begin{tabular}{ | c | c | }
\hline
  \multicolumn{2}{|c|}{$2^{nd}$-order Alpert Quadrature Rule for}   \\
  \multicolumn{2}{|c|}{integrals of the form $\int_{0}^{1}f(x)+g(x)\log(x)\, dx$,}   \\
  \multicolumn{2}{|c|}{with $a = 1$}\\
\hline
    NODES & WEIGHTS  \\
\hline
   1.591549430918953e-01  & 5.000000000000000e-01 \\
\hline
\end{tabular}
\end{center} \end{minipage}\\
\vspace{3em}
\begin{minipage}{0.48\linewidth} \begin{center}
\begin{tabular}{ | c | c | }
\hline
  \multicolumn{2}{|c|}{$6^{th}$-order Alpert Quadrature Rule for}   \\
  \multicolumn{2}{|c|}{integrals of the form $\int_{0}^{1}f(x)+g(x)\log(x)\, dx$,}   \\
  \multicolumn{2}{|c|}{with $a = 3$}\\
\hline
    NODES & WEIGHTS  \\
\hline
     4.004884194926570e-03  & 1.671879691147102e-02 \\
     7.745655373336686e-02  & 1.636958371447360e-01 \\
     3.972849993523248e-01  & 4.981856569770637e-01 \\
     1.075673352915104e+00  & 8.372266245578912e+00 \\
     2.003796927111872e+00  & 9.841730844088381e+00 \\
\hline
\end{tabular}
\end{center} \end{minipage}\\
\vspace{4em}
\begin{minipage}{0.48\linewidth} \begin{center}
\begin{tabular}{ | c | c | }
\hline
  \multicolumn{2}{|c|}{$10^{th}$-order Alpert Quadrature Rule for}   \\
  \multicolumn{2}{|c|}{integrals of the form $\int_{0}^{1}f(x)+g(x)\log(x)\, dx$,}   \\
  \multicolumn{2}{|c|}{with $a = 6$}\\
\hline
    NODES & WEIGHTS  \\
\hline
    1.175089381227308e-03 & 4.560746882084207e-03  \\
    1.877034129831289e-02 & 3.810606322384757e-02  \\
    9.686468391426860e-02 & 1.293864997289512e-01  \\
    3.004818668002884e-01 & 2.884360381408835e-01  \\
    6.901331557173356e-01 & 4.958111914344961e-01   \\
    1.293695738083659e+00 & 7.077154600594529e-01  \\
    2.090187729798780e+00 & 8.741924365285083e-01  \\
    3.016719313149212e+00 & 9.661361986515218e-01  \\
    4.001369747872486e+00 & 9.957887866078700e-01  \\
    5.000025661793423e+00 & 9.998665787423845e-01  \\
\hline
\end{tabular}
\end{center} \end{minipage}
\end{center}


\end{document}